%% file: smms_tractor.tex
\documentclass{amsart}
\usepackage{amsmath}
\usepackage{amssymb}
\usepackage{amsthm}

\DeclareMathOperator{\id}{id}

\DeclareMathOperator{\im}{im}
\DeclareMathOperator{\tr}{tr}

\DeclareMathOperator{\dvol}{dvol}

\DeclareMathOperator{\Ric}{Ric}

\DeclareMathOperator{\Rm}{Rm}

\DeclareMathOperator{\End}{End}
\DeclareMathOperator{\Hom}{Hom}

\DeclareMathOperator{\Ad}{Ad}

\DeclareMathOperator{\Hol}{Hol}
\DeclareMathOperator{\GL}{GL}

\DeclareMathOperator{\SO}{SO}

\newcommand{\og}{\overline{g}}

\newcommand{\lp}{\langle}
\newcommand{\rp}{\rangle}
\newcommand{\lv}{\lvert}
\newcommand{\rv}{\rvert}

\newcommand{\semiplus}{
  \mbox{$
  \begin{picture}(12.7,8)(-.5,-1)
  \put(2,0.2){$+$}
  \put(6.2,2.8){\oval(8,8)[l]}
  \end{picture}$}}
\newcommand{\contr}{\lrcorner}

\newcommand{\charconstant}{characteristic constant}

\newcommand{\mA}{\mathcal{A}}
\newcommand{\mB}{\mathcal{B}}
\newcommand{\mC}{\mathcal{C}}

\newcommand{\mE}{\mathcal{E}}

\newcommand{\mH}{\mathcal{H}}

\newcommand{\mL}{\mathcal{L}}
\newcommand{\mM}{\mathcal{M}}

\newcommand{\mP}{\mathcal{P}}
\newcommand{\mQ}{\mathcal{Q}}
\newcommand{\mR}{\mathcal{R}}
\newcommand{\mS}{\mathcal{S}}
\newcommand{\mT}{\mathcal{T}}

\newcommand{\kg}{\mathfrak{g}}

\newcommand{\bD}{\mathbb{D}}

\newcommand{\bN}{\mathbb{N}}

\newcommand{\bR}{\mathbb{R}}

\newcommand{\bT}{\mathbb{T}}

\newcommand{\J}{\mathrm{J}}
\newcommand{\bg}{\mathbf{g}}

\newtheorem{thm}{Theorem}[section]
\newtheorem{prop}[thm]{Proposition}
\newtheorem{lem}[thm]{Lemma}
\newtheorem{cor}[thm]{Corollary}

\theoremstyle{definition}
\newtheorem{defn}[thm]{Definition}

\newtheorem{example}[thm]{Example}
\newtheorem*{conv}{Convention}

\theoremstyle{remark}
\newtheorem{remark}[thm]{Remark}

\numberwithin{equation}{section}

\begin{document}

\title[SMMS, Quasi-Einstein Metrics, and Tractors]{Smooth metric measure spaces, quasi-Einstein metrics, and tractors}
\author{Jeffrey S. Case}
\thanks{Partially supported by NSF-DMS Grant No.\ 1004394}
\address{Department of Mathematics \\ Princeton University \\ Princeton, NJ 08544}
\email{jscase@math.princeton.edu}
%\date{}
\begin{abstract}
We introduce the tractor formalism from conformal geometry to the study of smooth metric measure spaces.  In particular, this gives rise to a correspondence between quasi-Einstein metrics and parallel sections of certain tractor bundles.  We use this formulation to give a sharp upper bound on the dimension of the vector space of quasi-Einstein metrics, providing a different perspective on some recent results of He, Petersen and Wylie.
\end{abstract}
\maketitle

%%%%%%%%%%%%%%%%%%%%%%%%%%%%%%%%%%%%%%%%%%%%%%%%%%%%%%%%%%%%%
%                                                           %
% Structure of the document                                 %
%                                                           %
% 1. Intro                                                  %
%   *. Acknowledgments                                      %
% 2. Tractor Bundles                                        %
% 3. Smooth metric measure spaces                           %
% 4. Quasi-Einstein Metrics                                 %
%   a. The Old Formulation                                  %
%   b. The New Formulation                                  %
%   c. The Linear Formulation                               %
% 5. Properties of $C$ and $W$-bundles                      %
%   a. Preliminary Computations                             %
%   b. Properties                                           %
%   c. Consequences for Parallel Sections                   %
% 6. $C$ and $W$-Holonomy                                   %
%   a. Framework                                            %
%   b. Results                                              %
% 7. Conclusion                                             %
% A. Appendix: R^n as a model                               %
% B. Appendix: Algebraic Formulation                        %
%                                                           %
%%%%%%%%%%%%%%%%%%%%%%%%%%%%%%%%%%%%%%%%%%%%%%%%%%%%%%%%%%%%%

\input{intro}

\subsection*{Acknowledgments}
I would like to first thank Robert Bartnik, who posed to me the question of whether the tractor calculus can be used to describe static metrics.  I would also like to thank Rod Gover for explaining to me many of the ideas underlying his almost Einstein metrics.  Additionally, I would like to thank Chenxu He, Peter Petersen, and especially William Wylie for many discussions of ``singular'' quasi-Einstein metrics and their recent work~\cite{HePetersenWylie2010,HePetersenWylie2011c}.  Finally, I would like to thank the referees for their many useful suggestions for improving this article.

\input{tractors}
\input{smms}
\input{formulation}
\input{bundle_properties}
\input{holonomy}
\input{conclusion}

\appendix
\input{model}
\input{algebra}

\bibliographystyle{abbrv}
\bibliography{../bib}
\end{document}

%% file: intro.tex
\section{Introduction}
\label{sec:intro}

Smooth metric measure spaces are Riemannian manifolds $(M^n,g)$ equipped with a smooth measure $e^{-\phi}\dvol_g$.  Given a parameter $m\in\bR\cup\{\pm\infty\}$, one defines the ($m$-)Bakry-\'Emery Ricci tensor
\[ \Ric_\phi^m := \Ric + \nabla^2\phi - \frac{1}{m}d\phi\otimes d\phi, \]
where $\Ric$ is the Ricci curvature of $g$ and $\nabla^2$ is the Hessian operator.  This tensor was introduced in the case $m\geq0$ by Bakry and \'Emery~\cite{BakryEmery1985} due to its role as the curvature term in the Bochner inequality on these spaces, and as such is the natural analogue of the Ricci tensor on smooth metric measure spaces.

Quasi-Einstein metrics are Riemannian metrics $g$ such that there exist $\phi\in C^\infty(M)$ and constants $\lambda,m$ such that $\Ric_\phi^m=\lambda g$.  Viewed this way, they are natural generalizations of Einstein metrics, and they include as special cases gradient Ricci solitons, the bases of warped product Einstein metrics, and the bases of conformally Einstein product manifolds, depending on the value of $m$ (see~\cite{Case2010a} and references therein).  For this reason, it is natural to regard $m$ as an interpolating constant, and to ask to what extent ideas used to study Einstein metrics or gradient Ricci solitons, say, can be generalized to study all quasi-Einstein metrics.

In~\cite{Case2010a}, the author described a natural way to define pointwise conformal transformations of smooth metric measure spaces.  This raises the possibility of using ideas from conformal geometry to better understand smooth metric measure spaces.  One success of this idea is found in~\cite{Case2010b}, where it was shown that the Yamabe constant and Perelman's entropies are special cases of a family of ``conformal invariants'' parameterized by $m$ which are associated to a smooth metric measure space.

The purpose of this article is to further develop the notion of a smooth metric measure space as a conformal object.  Our principal goal is to introduce the tractor calculus to this setting, following in many ways the foundational paper of Bailey, Eastwood and Gover~\cite{Bailey1994}.  More precisely, we will introduce a codimension one subbundle $\bT^W$ of the standard tractor bundle $\bT$ together with a new connection $\nabla^W$ and a new tractor-$D$ operator defined on $\bT$ which characterize conformally quasi-Einstein metrics as $\nabla^W$-parallel sections of $\bT^W$.  This perspective generalizes recent ideas of He, Petersen and Wylie~\cite{HePetersenWylie2010,HePetersenWylie2011c} for studying the vector space of quasi-Einstein metrics.  Moreover, by considering the singular sets of these tractors, we arrive at a new way to regard conformal infinities of Poincar\'e-Einstein manifolds and horizons of static metrics as ``the same'' (cf.\ \cite{Corvino2000,Gover2008}).  Additionally, the curvature of the connection $\nabla^W$ will give us natural analogues of the Weyl and the Cotton tensors for a smooth metric measure space, giving new meaning to the algebraic curvature tensors introduced in~\cite{Case2010b,HePetersenWylie2010} through \emph{ad hoc} means.

The main accomplishment of this article is the aforementioned correspondence between conformally quasi-Einstein metrics and $\nabla^W$-parallel sections of $\bT^W$.  Stated precisely, with terminology to be introduced in the ensuing sections, this correspondence is given as follows:

\begin{thm}
\label{thm:equivalence}
Let $(M^n,c,v^md\nu)$ be a smooth conformal measure space with \charconstant\ $\mu$ such that $n\geq3$ and $m\in\bR\setminus\{-n,1-n,2-n\}$.  There exists a codimension one subbundle $\bT^W\subset\bT$ of the standard tractor bundle, together with a connection $\nabla^W\colon\mT\to T^\ast M\otimes\mT$ and an operator $\bD^W\colon\mE[1]\to\mT$ with the following properties:
\begin{enumerate}
\item $u\in\mE[1]$ determines a quasi-Einstein scale wherever it is nonvanishing if and only if $\bD^Wu\in\mT$ is parallel with respect to $\nabla^W$.
\item If $I\in\mT^W$ is parallel with respect to $\nabla^W$, then there is a $u\in\mE[1]$ such that $I=\frac{1}{m+n}\bD^Wu$.
\end{enumerate}
\end{thm}

Our notation will be explained in Section~\ref{sec:tractors}, Section~\ref{sec:smms}, and Section~\ref{sec:formulation}.  In particular, Section~\ref{sec:tractors} provides a brief summary of the usual tractor calculus in conformal geometry, while Section~\ref{sec:smms} provides a brief summary of smooth metric measure spaces and introduces their analogues as objects in conformal geometry, which we term smooth conformal measure spaces.  These definitions are then brought together in Section~\ref{sec:formulation} to define the subbundle $\bT^W$ and the operators $\nabla^W$ and $\bD^W$.  The verification that these objects are all well-defined, as well as the proof of Theorem~\ref{thm:equivalence}, are then presented in Section~\ref{sec:bundle_properties}.

More accessible at the moment are the geometric interpretations of Theorem~\ref{thm:equivalence}, and in particular the assumption that the smooth conformal measure space $(M^n,c,v^md\nu)$ with \charconstant\ $\mu$ admits a quasi-Einstein scale $u\in\mE[1]$.  Fix a metric $g\in c$, which allows us to regard $u$ and $v$ as functions.  Then
\begin{enumerate}
\item If $m\in\bN$ and $(F^m,h)$ is any Einstein manifold with $\Ric_h=\mu h$, the metric $u^{-2}(g\oplus v^2h)$ on $M\times F$ is an Einstein metric with scalar curvature determined by the length of the corresponding tractor $I\in\mT^W$.
\item If $p=2-m-n\in\bN$ and $(F^p,h)$ is an Einstein manifold with scalar curvature determined by the length of $I$, the metric $v^{-2}(g\oplus u^2h)$ on $M\times F$ is an Einstein metric with scalar curvature $(n+p)\mu$.
\end{enumerate}
In particular, when $m>0$, we see that scale singularities (zero sets) of $u$ represent conformal infinities, while when $m<2-n$, the scale singularities of $u$ are the sets where the warping function degenerates, which are known as horizons in the case $m=1-n$ (cf.\ \cite{Corvino2000}); this is the aforementioned relationship between conformal infinities and horizons.  Also, the above correspondence reveals the limitations of Theorem~\ref{thm:equivalence}: Besides being unable to treat gradient Ricci solitons (corresponding to $m=\pm\infty$), Theorem~\ref{thm:equivalence} does not yield a tractor formulation of the problem of finding the warping function for an Einstein warped product with a fiber of dimension one or two.  The difficulties in these cases arise because certain formulae we give become singular at these values, and it is not clear that these singularities can be removed.

Since Theorem~\ref{thm:equivalence} yields a correspondence between quasi-Einstein metrics and tractors which are parallel with respect to a certain connection, it is unsurprising that we can use the holonomy groups of these connections to better understand the structure of manifolds which admit quasi-Einstein metrics.  In particular, we can establish sharp bounds on the number of linearly independent quasi-Einstein scales on a smooth conformal measure space, analogous to a result of He, Petersen and Wylie~\cite{HePetersenWylie2011c}; for a more precise comparison of our results to those given in~\cite{HePetersenWylie2011c}, see Remark~\ref{rk:holonomy}.

The tractor calculus has many important applications in conformal geometry, and for this reason, we likewise expect that it will have important applications to the study of smooth metric measure spaces in general, and quasi-Einstein metrics in particular.  Some evidence for this can be found in the article~\cite{Case2011o}, where the construction by Gover and Nurowski~\cite{GoverNurowski2006} of obstructions to the existence of almost Einstein metrics in a conformal class is generalized to the setting of smooth metric measure spaces, yielding analogous obstructions to the existence of quasi-Einstein scales.  Some additional insights into smooth metric measure spaces granted by our tractor formulation, as well as a brief discussion of an open problem in the study of quasi-Einstein metrics which is naturally related to our perspective, are discussed in Section~\ref{sec:conclusion}.

This article is organized as follows.

In Section~\ref{sec:tractors}, we introduce those aspects of the tractor calculus which we shall use here.  We elect to work directly with vector bundles as much as possible, mostly following the presentation given in~\cite{Bailey1994}.

In Section~\ref{sec:smms}, we recall some basic definitions and ideas involved in the study of smooth metric measure spaces and quasi-Einstein metrics, and make precise the notion of a smooth conformal measure space.

In Section~\ref{sec:formulation}, we describe two ways to discuss quasi-Einstein metrics in the tractor language.  First, in Section~\ref{sec:formulation/new}, we describe how this can be done without the notion of a smooth conformal measure space.  Second, in Section~\ref{sec:formulation/linear}, we specialize to the setting of smooth conformal measure spaces and describe the bundle and operators of Theorem~\ref{thm:equivalence}.

In Section~\ref{sec:bundle_properties}, we establish some basic properties of the subbundle $\bT^W$ and the operators $\nabla^W$ and $\bD^W$.  In particular, we show that $\nabla^W$ and $\bD^W$ are well-defined tractor operators and give the proof of Theorem~\ref{thm:equivalence}.  Additionally, we observe that the singular sets of sections of $\bT^W$ with respect to $\nabla^W$ are well-behaved in the expected way (cf.\ \cite{Gover2008}).  In particular, this gives a simple way to understand the geometric similarity of horizons and conformal infinities.

In Section~\ref{sec:holonomy}, we consider the holonomy of $(\bT,\nabla^W)$, and use it to give new proofs of some recent results of He, Petersen and Wylie~\cite{HePetersenWylie2011c} for Riemannian manifolds admitting ``many'' quasi-Einstein metrics.

In Section~\ref{sec:conclusion}, we conclude with discussions relating our work to other work in conformal geometry and smooth metric measure spaces.  In particular, we give a very brief summary of~\cite{HePetersenWylie2011c} and ask whether holonomy methods yield new insights into that work and the related question of finding the global classification of quasi-Einstein metrics.  We also observe that our tractor formulation yields immediately a natural analogue of the Paneitz operator on smooth metric measure spaces, and note that our work also can be used to describe a class of critical metrics recently studied by Miao and Tam~\cite{MiaoTam2008,MiaoTam2011}.

We also include two appendices to explain a few remarks made throughout this article.  In Appendix~\ref{sec:model}, we recall the standard tractor bundle and its natural basis of parallel sections on the standard $n$-sphere, and use it to give another perspective on some important examples of quasi-Einstein metrics.  In Appendix~\ref{sec:algebra}, we briefly discuss some of the underlying algebraic structure of the tractor calculus, and use it to describe in what sense the connection $\nabla^W$ is canonical.

%% file: tractors.tex
\section{Tractor Bundles}
\label{sec:tractors}

We begin by recalling some of the basic ingredients of the tractor calculus in conformal geometry.  In order to keep our discussion simple and to emphasize the computational utility of the tractor bundle, we shall restrict our discussion to the point of view afforded by vector bundles, following Bailey, Eastwood and Gover~\cite{Bailey1994}.  For an alternative treatment in the more general framework of parabolic geometries, we refer the interested reader to~\cite{CapSlovak2009}.

Just as the basic objects of Riemannian geometry are smooth functions and tensor fields, the basic objects of conformal geometry are densities and tractors.  The fundamental object of conformal geometry is a conformal class $c$ of metrics,
\[ c = [g] = \{ h \colon h=e^{2s}g \mbox{ for some } s\in C^\infty(M) \} \]
for some (and hence any) $g\in c$.  Throughout this article we shall only focus on the case that $c$ is a Riemannian conformal class; i.e.\ elements $g\in c$ are Riemannian metrics.

Alternatively, one can regard $c$ as specifying a ray subbundle $\mC\subset S^2T^\ast M$ by requiring that sections of $\mC$ are metrics in $c$.  Regarded this way, the natural actions $\delta_s\colon\mC\to\mC$ given by $\delta_s(x,g_x)=(x,s^2g_x)$ for any $s>0$ give $\mC$ the structure of a $\bR^+$-principle bundle.  This allows one to introduce conformal density bundles.

\begin{defn}
Let $(M,c)$ be a conformal manifold and fix $w\in\bR$.  The \emph{(conformal) density bundle of weight $w$} is the line bundle $E[w]$ associated to $\mC$ via the representation $\bR^+\ni t\mapsto t^{-w/2}\in\End(\bR)$.

We denote by $\mE[w]=\Gamma\left(E[w]\right)$, the space of smooth sections of $E[w]$.
\end{defn}

The conformal density bundles are trivial line bundles, and \emph{choice of scale} $g\in c$ defines an isomorphism $\mE[w]\cong C^\infty(M)$, denoted $\sigma\cong_g\sigma_g$, by $\sigma_g(x)=\sigma(x,g_x)$.  Given another choice of scale $h=e^{2s}g$, the functions $\sigma_g$ and $\sigma_h$ are related by $\sigma_h=e^{ws}\sigma_g$.  We shall henceforth denote this relationship by $\sigma\mapsto e^{ws}\sigma$.

Given a vector bundle $V$ over $M$, one can consider the tensor product $V\otimes E[w]$, which we shall denote by $V[w]$.  In this way, one easily sees that the conformal class $c$ determines the \emph{conformal metric} $\bg\in S^2T^\ast M[2]$ by $\bg(\cdot,g):=g$.  This section determines the canonical isomorphism $TM[w]\cong T^\ast M[w+2]$, which we shall use without comment.

The next important object in the tractor calculus is the standard tractor bundle, which in many ways serves as the analogue of the tangent bundle in Riemannian geometry.

\begin{defn}
\label{defn:tractor}
Let $(M^n,c)$ be a conformal manifold.  The \emph{standard tractor bundle} $\bT$ is the vector bundle $M$ characterized by the following properties.
\begin{enumerate}
\item Any choice of metric $g\in c$ induces an isomorphism $\bT\cong \bR\oplus TM\oplus \bR$, where $\bR$ denotes the trivial (real) line bundle over $M$.  Given $I\in\bT$, we shall denote this isomorphism by $I=(\rho,\omega,\sigma)$ for $\sigma,\rho\in\bR$ and $\omega\in TM$, or equivalently,
\begin{equation}
\label{eqn:tractor_notation}
I = \begin{pmatrix}\sigma\\\omega\\\rho\end{pmatrix} .
\end{equation}
\item The isomorphism $\bT\cong\bR\oplus TM\oplus\bR$ transforms with a change of scale as
\begin{equation}
\label{eqn:tractor_law}
\begin{pmatrix} \sigma\\\omega\\\rho \end{pmatrix} \mapsto \begin{pmatrix} e^s\sigma\\e^{-s}(\omega+\sigma\,\nabla s)\\e^{-s}(\rho-g(\nabla s,\omega) - \frac{1}{2}|\nabla s|^2\,\sigma) \end{pmatrix} ,
\end{equation}
where $\nabla s$ and $|\nabla s|^2$ are computed with respect to $g$.
\end{enumerate}

A \emph{(standard) tractor} is a section of $\mT$.
\end{defn}

\begin{conv}
Unless otherwise specified, whenever we use the symbol $I$ to denote a tractor and a choice of scale $g\in c$ is understood, we will use the symbols $\sigma,\omega,\rho$ to denote its components as in~\eqref{eqn:tractor_notation}.
\end{conv}

\begin{remark}
\label{rk:plus_notation}
In the above definition, the transformation formula~\eqref{eqn:tractor_law} implies that a choice of scale $g\in c$ induces a splitting $\bT\cong E[1]\oplus TM[-1]\oplus E[-1]$ of vector bundles, but that this splitting is not canonical, in that it cannot be made to depend only on $c$.  However, reading from the top down, the first nonzero entry of the above vector representation of a tractor $I$ is conformally invariant, and will be referred to as the \emph{projecting part}.  Structurally, we also see that the projection $\bT\to E[1]$ is well-defined, as is the inclusion $E[-1]\hookrightarrow\bT$.  Following~\cite{BransonGover2005}, we will record these observations by writing $\mT=\mE[1]\semiplus TM[-1]\semiplus\mE[-1]$.
\end{remark}

\begin{remark}
One can take tensor products with $\bT$ in the obvious way.  An important example is
\[ T^\ast M\otimes\bT = T^\ast M[1] \semiplus (T^\ast M\otimes TM[-1]) \semiplus T^\ast M[-1], \]
which arises when discussing connections on $\bT$.  Using the obvious ``coordinate'' representation of a section of $T^\ast M\otimes\bT$ given a choice of scale, one sees that the transformation law~\eqref{eqn:tractor_law} implies that
\[ \begin{pmatrix} \alpha\\T\\\beta \end{pmatrix} \mapsto \begin{pmatrix} e^s\alpha\\e^{-s}(T+\alpha\otimes\nabla s)\\e^{-s}(\beta-T(\cdot,ds)-\frac{1}{2}|\nabla s|^2\alpha) \end{pmatrix}, \]
where $\alpha,\beta\in T^\ast M$ and $T\in T^\ast M\otimes TM$.
\end{remark}

Given $(M^n,c)$ the standard tractor bundle comes equipped with a canonical metric and connection.  Here and in the following, we use the same symbol to denote a tensor bundle and its space of smooth sections; for example, $T^\ast M$ will denote either the cotangent bundle or the space of one-forms, depending on the context.  Also, given two vector bundles $V$ and $W$, we write $\Gamma(V)\otimes\Gamma(W):=\Gamma(V\otimes W)$; in particular, $\mT^\ast\otimes\mT^\ast$ is the space of smooth sections of $\bT^\ast\otimes\bT^\ast$.

\begin{defn}
Let $(M^n,c)$ be a conformal manifold and let $\bT$ be the standard tractor bundle.

The \emph{tractor metric} $h\in\mT^\ast\otimes\mT^\ast$ is defined in the scale $g\in c$ by
\[ h(I,I) = 2\sigma\rho + |\omega|_g^2 . \]

The \emph{normal tractor connection} $\nabla\colon\mT\to T^\ast M\otimes\mT$ is defined in the scale $g\in c$ by
\[ \nabla I = \begin{pmatrix} \nabla\sigma-g(\omega,\cdot)\\\nabla\omega + \sigma P + \rho\,g\\\nabla \rho - P(\omega,\cdot) \end{pmatrix} , \]
where on the right-hand side, $\nabla$ is the Levi-Civita connection of $g$, $P=\frac{1}{n-2}\left(\Ric - \J\,g\right)$ is the \emph{Schouten tensor} with $\J=\tr P$ its trace, and in the second line we use the isomorphism $T^\ast M\otimes T^\ast M\cong T^\ast M\otimes TM[-2]$.
\end{defn}

Straightforward computations yield that the tractor metric and the normal tractor connection are both well-defined; i.e.\ they both transform under a change of scale as required by~\eqref{eqn:tractor_law}.  It is also clear that $h$ has signature $(n+1,1)$.  Since the normal tractor connection differs (in scale) from the Levi-Civita connection by a tensor, it is clear that it is indeed a connection.  Indeed, the normal tractor connection is effectively the normal Cartan connection (cf.\ \cite{CapSlovak2009}).  From this, or by direct computation, one has that the normal tractor connection is a metric connection, $\nabla h=0$.  Hence we may freely identify $\bT\cong\bT^\ast$, as we shall do without further comment.

\begin{remark}
Our conventions in writing $I$ in scale in Definition~\ref{defn:tractor} are made so that $h(I,J)$ can equivalently be computed by taking the usual matrix product of $I$ as a row vector and $J$ as a column vector.  We will also find it expedient to use the notation $\lp I,J\rp$ and $\lv I\rv^2$ to denote the quantities $h(I,J)$ and $h(I,I)$, respectively.  Especially in the latter case, one must remember that $h$ has signature $(n+1,1)$, so that $\lv I\rv^2$ is not in general nonnegative.
\end{remark}

Using the tractor metric, our earlier observation that the projection $\bT\to E[1]$ and the inclusion $E[-1]\hookrightarrow\bT$ are well-defined gives rise to an important tractor.

\begin{defn}
Let $(M,c)$ be a conformal manifold.  Given a choice of scale $g\in c$, the \emph{projector} $X\in\mT[1]$ is the tractor $X=(1,0,0)$.
\end{defn}

As a consequence of~\eqref{eqn:tractor_law}, it is easily seen that $X$ is well-defined.  The projection $\bT\to E[1]$ is given by $I\mapsto h(I,X)$, while the inclusion $E[-1]\hookrightarrow\bT$ is given by $\rho\mapsto X\otimes\rho$.

Another important operator is the tractor-$D$ operator.

\begin{defn}
Let $(M^n,c)$ be a conformal manifold.  The \emph{tractor-$D$ operator} $\bD\colon\mE[w]\to\mT[w-1]$ is defined by
\[ \bD\sigma = \begin{pmatrix} w(n+2w-2)u\\(n+2w-2)\nabla u\\-\left(\Delta u+w\J u\right)\end{pmatrix} \]
in any choice of scale $g\in c$.
\end{defn}

A straightforward computation shows that the tractor-$D$ operator is well-defined; i.e.\ it transforms with a change of scale according to~\eqref{eqn:tractor_law}.

\begin{remark}
The tractor-$D$ operator easily extends to a map between arbitrary tractor bundles by replacing the connection and the Laplacian with the corresponding operators induced by the normal tractor connection.  In this level of generality, the tractor-$D$ operator plays an important role in describing conformally covariant operators; see~\cite{GoverPeterson2003}.
\end{remark}

We will also find it useful to introduce the adjoint tractor bundle, which for our purposes will be useful for discussing other connections on $\mT$ and the curvature of those connections.

\begin{defn}
\label{defn:adjoint}
Let $(M^n,c)$ be a conformal Riemannian manifold.  The \emph{adjoint tractor bundle} is the bundle $\mA:=\bT\wedge\bT$, together with
\begin{enumerate}
\item the contraction map $\mA\otimes\mT\mapsto\mT$, denoted $A(I):=I\contr A$, with the convention
\[ K\contr (I\wedge J) = \lp I,K\rp J - \lp J,K\rp I, \]
\item the metric $h\in\mA^\ast\otimes\mA^\ast$ induced from the tractor metric
\[ h(I_1\wedge I_2,J_1\wedge J_2) = \lp I_1,J_1\rp \lp I_2,J_2\rp - \lp I_1,J_2\rp \lp I_2,J_1\rp , \qquad \mbox{and} \]
\item the Lie algebra bracket $\{\,,\,\}\colon\mA\times\mA\to\mA$
\[ \{A,B\}(I,J) = \lp A(I),B(J)\rp - \lp A(J),B(I)\rp . \]
\end{enumerate}
\end{defn}

Note that, as opposed to the standard tractor bundle $\bT$, we will use the symbol $\mA$ to denote both the adjoint tractor bundle and its space of smooth sections.  We will maintain the distinction between the standard tractor bundle $\bT$ and its space of smooth sections $\mT$ for clarity in discussing holonomy, an issue we will not have when discussing other vector bundles.

Using the facts $\bT\wedge\bT=TM[0]\semiplus(\Lambda^2T^\ast M[2]\oplus\mE[0])\semiplus T^\ast M[0]$ (cf.\ \cite{Bailey1994}) and $\Lambda^2 TM\oplus\bR\cong \mathfrak{co}(n)$, it is straightforward to check that this is (isomorphic to) the usual adjoint tractor bundle introduced by \v{C}ap and Gover~\cite{CapGover2002} (see also Appendix~\ref{sec:algebra}).  Moreover, just as the important tractor $X\in\mT[1]$ appeared as a consequence of the decomposition of $\mT$, there is an important pseudo-tractor which appears as a consequence of the decomposition of $\mA$.

\begin{defn}
Let $(M^n,c)$ be a conformal manifold and choose a scale $g\in c$.  The \emph{$Z$-projector} $Z_g\in T^\ast M\otimes\mT_g[1]$ is given by $Z_g(\omega)=(0,\omega,0)$ for all $\omega\in TM[0]$.
\end{defn}

While $Z_g$ is not actually a tractor, as it does not satisfy the conformal transformation law~\eqref{eqn:tractor_law}, it is important for two reasons.  First, given a choice of scale, $Z_g=\nabla X$.  Second, $Z_g$ can be used to realize the projection $\mA\to TM[0]$ and the inclusion $T^\ast M[0]\hookrightarrow\mA$.

\begin{prop}
Let $(M^n,c)$ be a conformal manifold with adjoint tractor bundle $(\mA,h)$.  Then the tractor $X\wedge Z\in T^\ast M\otimes\mA[2]$, given by $X\wedge Z_g$ for any choice of scale, is such that
\begin{enumerate}
\item the map
\[ \mA \ni A \mapsto h(A,X\wedge Z) \in T^\ast M[2] \cong TM[0] \]
is the projection map; and
\item the map
\[ T^\ast M[0] \ni \omega \mapsto X\wedge Z(\omega) \in \mA \]
is the inclusion map, where we regard $\omega\in TM[-2]$ on the right hand side via $T^\ast M[0]\cong TM[-2]$.
\end{enumerate}
\end{prop}

\begin{proof}

Since $X\wedge X=0\in\mA[2]$ and $Z_h=Z_g-ds\otimes X$, it follows immediately that $X\wedge Z\in T^\ast M\otimes\mA[2]$ is well-defined.  The final two claims follow directly from the definitions of $X$ and $Z$.
\end{proof}

Using the inclusion $T^\ast M\hookrightarrow\mA$, the contraction $\mA\otimes\mT\mapsto\mT$ restricts to an action $T^\ast M\otimes\mT\mapsto\mT$, given by
\[ \alpha\otimes\begin{pmatrix}u\\\omega\\\rho\end{pmatrix} \mapsto \begin{pmatrix}0\\u\alpha^\sharp\\-\alpha(\omega)\end{pmatrix}, \]
where $\sharp\colon T^\ast M\to TM[-2]$ is the isomorphism induced by $\bg$.  In this way, we see that this action is precisely Kostant's codifferential $\partial^\ast\colon T^\ast M\otimes\mT\to\mT$ (cf.\ \cite{CapSlovak2009}).  This action plays an important role in characterizing connections on tractor bundles; see Appendix~\ref{sec:algebra} for further details.

%% file: smms.tex
\section{Smooth Metric Measure Spaces}
\label{sec:smms}

Let us briefly recall some basic definitions of smooth metric measure spaces and natural geometric definitions associated to these objects.  This will culminate in our definition of smooth conformal measure spaces.

\begin{defn}
\label{defn:smms}
A \emph{smooth metric measure space (SMMS)} is a four-tuple $(M^n,g,v^m\dvol_g,m)$ of a smooth Riemannian manifold $(M^n,g)$ with its Riemannian volume element $\dvol_g$, a positive function $v\in C^\infty(M)$ determining a \emph{smooth measure} $v^m\dvol_g$, and a dimensional parameter $m\in\bR\cup\{\pm\infty\}$.
\end{defn}

Two comments about this definition are in order.  First, we adopt the convention that $v^0=1$, so that SMMS with $m=0$ are simply Riemannian manifolds; i.e.\ we are not imposing additional data in this case.  Second, the above definition does not make sense when $\lv m\rv=\infty$ as stated.  We overcome this by associating to an SMMS the function $\phi$ (formally) defined by $v^m=e^{-\phi}$, so that an SMMS can equivalently be regarded as a four-tuple $(M^n,g,e^{-\phi}\dvol_g,m)$ with the assumption $\phi\equiv0$ if $m=0$, a definition which does make sense when $\lv m\rv=\infty$.

When writing an SMMS in an abstract sense, we will usually denote it by the triple $(M^n,g,v^m\dvol_g)$, with the dimensional parameter $m$ being encoded in our notation for the measure.  We will also always associate to an SMMS the function $\phi$ defined formally by $v^m=e^{-\phi}$ as above.  In this way, it will be easy to see the validity of definitions of certain geometric objects defined on an SMMS $(M^n,g,v^m\dvol)$ in the limiting cases $\lv m\rv=\infty$ (cf.\ Definition~\ref{defn:curvature}).

% When writing SMMS in an abstract sense, we will usually denote them by the triple $(M^n,g,v^m\dvol_g)$, with the dimensional parameter $m$ being encoded in our notation for the measure.  Note in particular that our conventions force the measure to be the standard Riemannian measure when $m=0$, so that we are not adding any new data in this case.  We define the function $\phi\in C^\infty(M)$ by $v^m=e^{-\phi}$; in this way we can easily make sense of the definition of an SMMS when $m=\pm\infty$.  Throughout this article, the symbol $\phi$ will be related to $v$ this way.

Roughly speaking, the role of the dimensional parameter $m$ is to specify that the measure $v^m\dvol_g$ should be regarded as an $(m+n)$-dimensional measure.  This is made precise through various geometric definitions involving SMMS, as we shall see below.  The motivation for this perspective is easily seen through the following example.

\begin{example}
\label{ex:cc}
Let $(M^n,g)$ and $(F^m,h)$ be smooth manifolds and fix a positive function $v\in C^\infty(M)$.  For any sequence $\{\varepsilon_i\}$ with $\varepsilon_i\to0$ as $i\to\infty$, the sequence of warped product manifolds
\begin{equation}
\label{eqn:cc_wp}
\left( M \times F, g\oplus (\varepsilon_i v)^2 h \right)
\end{equation}
converges in the sense of Cheeger-Colding~\cite{CheegerColding1997} to the SMMS $(M^n,g,v^m\dvol)$; i.e.\ the renormalized volume element of~\eqref{eqn:cc_wp} converges to the measure $v^m\dvol_g$ on $M$ as $i\to\infty$ (cf.\ Example~1.24 in~\cite{CheegerColding1997}).
\end{example}

In other words, one can always regard an SMMS $(M^n,g,v^m\dvol_g)$ with $m\in\bN$ as a collapsed limit of some $(m+n)$-dimensional warped product.  In this special case, all of our definitions can be realized by considerations involving the warped product $(M^n\times F^m,g\oplus v^2h)$.  However, it is important to note that our definitions both make sense and are useful for all values of $m$, including negative and infinite values, and so we will always state our definitions intrinsically; i.e.\ they will only ever depend on the SMMS $(M^n,g,v^m\dvol)$.

Our goal is to study natural ``weighted'' objects defined on an SMMS; i.e.\ objects which depend in a natural way on the specified metric and (weighted) measure.  An obvious such object is the weighted divergence.

\begin{defn}
Let $(M^n,g,v^m\dvol)$ be an SMMS.  The \emph{weighted divergence $\delta_\phi$} is the (negative of the) adjoint of the exterior derivative with respect to the measure $v^m\dvol$,
\[ \int_M \lp\alpha,d\beta\rp e^{-\phi}\dvol = -\int_M \lp\delta_\phi\alpha,\beta\rp e^{-\phi}\dvol \]
for all $\alpha\in\Lambda^{k+1}T^\ast M$, $\beta\in\Lambda^kT^\ast M$, with $\lp\cdot,\cdot\rp$ the natural (pointwise) inner product on $k$-forms induced by $g$.

The \emph{weighted Laplacian $\Delta_\phi\colon C^\infty(M)\to C^\infty(M)$} is defined by $\Delta_\phi=\delta_\phi d$.
\end{defn}

Relative to the Laplacian $\Delta=\tr_g\nabla^2$ defined by the Riemannian metric, the weighted Laplacian satisfies $\Delta_\phi u = \Delta u - \lp\nabla u,\nabla\phi\rp$ for all $u\in C^\infty(M)$.

Besides the weighted Laplacian $\Delta_\phi$, the central object when studying SMMS is the Bakry-\'Emery Ricci curvature, which was introduced as the natural curvature term in the Bochner formula for $\Delta_\phi\lv\nabla w\rv^2$; see~\cite{BakryEmery1985}.

\begin{defn}
\label{defn:curvature}
Let $(M^n,g,v^m\dvol)$ be an SMMS.  The \emph{Bakry-\'Emery Ricci tensor $\Ric_\phi^m$} is the tensor
\[ \Ric_\phi^m = \Ric - mv^{-1}\nabla^2v = \Ric + \nabla^2\phi - \frac{1}{m}d\phi\otimes d\phi . \]
The \emph{weighted scalar curvature $R_\phi^m$} is the tensor
\[ R_\phi^m = R - 2mv^{-1}\Delta v - m(m-1)v^{-2}\lv\nabla v\rv^2 = R + 2\Delta\phi - \frac{m+1}{m}\lv\nabla\phi\rv^2 . \]
\end{defn}

Note that, in general, $R_\phi^m\not=\tr_g\Ric_\phi^m$.  The relationship between the two is easily seen in the case $m\in\bN$ by considering the warped product $(M^n\times F^m,\og=g\oplus v^2h)$: The Bakry-\'Emery Ricci tensor is the horizontal part of $\Ric_{\og}$ (i.e.\ its restriction to lifts of vector fields on $M$), and if one assumes additionally that $h$ is scalar flat, then $R_\phi^m$ is the scalar curvature of $\og$; in other words, $R_\phi^m$ incorporates the extra information ``hidden'' in the Ricci curvature of the vertical directions.

The notion of a quasi-Einstein metric is the natural generalization of an Einstein metric made by using the Bakry-\'Emery Ricci tensor.

\begin{defn}
\label{defn:smms_qe}
An SMMS $(M^n,g,v^m\dvol)$ is said to be \emph{quasi-Einstein} if there exists a constant $\lambda\in\bR$ such that $\Ric_\phi^m=\lambda g$.  We call $\lambda$ the \emph{quasi-Einstein constant}.
\end{defn}

An important fact about quasi-Einstein metrics is that they determine a type of ``integrability condition.''

\begin{prop}[\cite{Kim_Kim}; see also~\cite{Case2010a}]
\label{prop:kk}
Let $(M^n,g,v^m\dvol)$ be a quasi-Einstein SMMS with quasi-Einstein constant $\lambda$.  Then there exists a constant $\mu\in\bR$ such that
\begin{equation}
\label{eqn:kk_csc}
R_\phi^m + m\mu v^{-2} = (m+n)\lambda .
\end{equation}
\end{prop}

Geometrically, the constant $\mu$ is such that whenever $(M^n,g,v^m\dvol)$ is a quasi-Einstein SMMS with $m\in\bN$ and $\mu$ as in~\eqref{eqn:kk_csc}, given \emph{any} Einstein manifold $(F^m,h)$ satisfying $\Ric_h=\mu h$, the warped product $(M^n\times F^m,\og=g\oplus v^2h)$ is Einstein with $\Ric_{\og}=\lambda\og$.  In other words, the quasi-Einstein condition implies that $\Ric_{\og}=\lambda\og$ when restricted to horizontal vector fields, while the condition~\eqref{eqn:kk_csc} implies that $\Ric_{\og}=\lambda\og$ on vertical vector fields (cf.\ \cite{Besse,Case_Shu_Wei,Kim_Kim}).  It is in this sense that we regard~\eqref{eqn:kk_csc} as an integrability condition; alternatively, one can regard it as an additional constraint on $v$, which is more in line with its use in~\cite{Kim_Kim}.

Because of the importance of the constant $\mu$, we will find it useful to give it a name.

\begin{defn}
Let $(M^n,g,v^m\dvol)$ be a quasi-Einstein SMMS with quasi-Einstein constant $\lambda$.  The \emph{\charconstant} is the constant $\mu$ determined by~\eqref{eqn:kk_csc}.
\end{defn}

One of the key insights of~\cite{Case2010a} is that there is a natural notion of a conformal transformation of an SMMS determined by the dimensional parameter $m$.

\begin{defn}
Two SMMS $(M^n,g,v^m\dvol_g)$ and $(M^n,\hat g,\hat v^m\dvol_{\hat g})$ are said to be \emph{(pointwise) conformally equivalent} if there is a positive function $u\in C^\infty(M)$ such that
\begin{equation}
\label{eqn:scms}
\left( M^n, \hat g,\hat v^m\dvol_{\hat g}\right) = \left( M^n, u^{-2}g, u^{-m-n}v^m\dvol_g \right) .
\end{equation}
\end{defn}

In other words, we insist that the measure $v^m\dvol_g$ transforms as would the Riemannian volume element of an $(m+n)$-dimensional manifold.  Note that one can equivalently write~\eqref{eqn:scms} as
\[ \left( M^n,\hat g,\hat v^m\dvol_{\hat g}\right) = \left( M^n, u^{-2}g, (u^{-1}v)^m\dvol_{u^{-2}g} \right) . \]
In this way, \eqref{eqn:scms} determines an equivalence relation on SMMS.  Equivalence classes for this relation are smooth conformal measure spaces.

\begin{defn}
\label{defn:scms}
A \emph{smooth conformal measure space (SCMS)} is a triple $(M^n,c,v^md\nu)$ of a conformal manifold $(M^n,c)$ together with the conformal volume element $d\nu\in\Lambda^nT^\ast M[n]$, a positive density $v\in\mE[1]$, and a dimensional parameter $m\in\bR$.
\end{defn}

The conformal volume element $d\nu\in\Lambda^nT^\ast M[n]$ is defined by $d\nu(\cdot,g)=\dvol_g$; i.e.\ it is the (conformally weighted) volume element of the conformal metric $\bg$.

In particular, note that when $m=0$, an SCMS is simply a conformal manifold: In this case, the measure $v^0d\nu=d\nu$ is simply the conformal volume element defined by $c$, and thus we are not imposing any additional structure on the conformal manifold $(M^n,c)$.

\begin{remark}
The conformal transformation formula~\eqref{eqn:scms} still makes sense when $m\in\{\pm\infty\}$ by (formally) defining the function $f$ by $u^{m+n-2}=e^f$.  In this case, a ``conformal transformation'' is simply a change of measure.  This is certainly a useful perspective on SMMS (cf.\ \cite{Case2010a,Case2010b}).  However, since the metric $g$ remains fixed under such a transformation, the equivalence classes determined by~\eqref{eqn:scms} depend on a metric, and not a conformal class, which is the reason we require that $m\in\bR$ in that definition.
\end{remark}

From the standpoint of conformal geometry, the real interest in Definition~\ref{defn:scms} is the dimensional parameter $m$:  It explicitly states that $v\in\mE[1]$ is \emph{not} to be regarded as fixing a metric $v^{-2}\bg\in c$, but rather as determining a conformal class $[c\oplus v^2h]$ of metrics on a product $M^n\times F^m$ where $h$ is some fixed (but arbitrary) Riemannian metric on $F$.

For the majority of this article, we will be interested in studying conformally invariant objects associated to an SMMS in the following sense.

\begin{defn}
Let $(M^n,g,v^m\dvol_g)$ be an SMMS.  A geometric invariant $T=T[g,v^m\dvol_g]$ defined in terms of $g$ and $v$ is said to be \emph{conformally invariant of weight $w$} if for any $s\in C^\infty(M)$,
\[ T\left[e^{2s}g,(e^sv)^m\dvol_{e^{2s}g}\right] = e^{ws}T[g,v^m\dvol_g] . \]
\end{defn}

In other words, if $T$ is a conformal invariant of an SMMS $(M^n,g,v^m\dvol_g)$, it can be regarded as a well-defined geometric invariant on the corresponding SCMS $(M^n,c,v^md\nu)$.

To pursue a tractor formulation of quasi-Einstein metrics, it is useful to know how the Bakry-\'Emery Ricci tensor and the weighted scalar curvature transform with a conformal change of metric.

\begin{lem}[\cite{Case2010a}]
\label{lem:conf_trans}
Let $(M^n,g,v^m\dvol_g)$ be an SMMS and let $u\in C^\infty(M)$ be a positive function.  The Bakry-\'Emery Ricci tensor $\Ric_{f,\phi}^m$ and the weighted scalar curvature $u^2R_{f,\phi}^m$ of the SMMS $(M^n,\hat g,\hat v^m\dvol_{\hat g})$ determined by~\eqref{eqn:scms} are given by
\begin{align*}
\Ric_{f,\phi}^m & = \Ric_\phi^m + (m+n-2)u^{-1}\nabla^2 u + \left(u^{-1}\Delta_\phi u - (m+n-1)u^{-2}\lv\nabla u\rv^2\right)g \\
R_{f,\phi}^m & = R_\phi^m + 2(m+n-1)u^{-1}\Delta_\phi u - (m+n)(m+n-1)u^{-2}\lv\nabla u\rv^2,
\end{align*}
respectively, where all derivatives and traces on the right hand side are computed with respect to $g$.
\end{lem}

\begin{remark}
The notation $\Ric_{f,\phi}^m$ and $R_{f,\phi}^m$ is used to indicate that these formulae are well-defined in the limits $\lv m\rv=\infty$ by formally defining $u^{m+n-2}=e^f$; see~\cite{Case2010a} for details.
\end{remark}

In particular, one can easily derive the formulae specifying that an SMMS be conformally quasi-Einstein.

\begin{prop}[\cite{Case2010a}]
\label{prop:conf_qe}
Let $(M^n,g,v^m\dvol_g)$ be an SMMS and suppose that $u\in C^\infty(M)$ is such that the SMMS $(M^n,\hat g,\hat v^m\dvol_{\hat g})$ determined by~\eqref{eqn:scms} is a quasi-Einstein SMMS with quasi-Einstein constant $\lambda$ and \charconstant\ $\mu$.  Then
\begin{subequations}
\label{eqn:qe_conf}
\begin{align}
\label{eqn:qe_tf} 0 & = \left(uv\Ric + (m+n-2)v\nabla^2 u - mu\nabla^2 v\right)_0 \\
\label{eqn:qe_lambda} n\lambda v^2 & = (uv)^2R + (m+2n-2)uv^2\Delta u - mu^2v\Delta v \\
\notag & \quad - (m+n-1)nv^2\lv\nabla u\rv^2 + mnuv\lp\nabla u,\nabla v\rp \\
\label{eqn:qe_mu} n\mu u^2 & = (uv)^2R + (m+n-2)uv^2\Delta u - (m-n)u^2v\Delta v \\
\notag & \quad - (m+n-2)nuv\lp\nabla u,\nabla v\rp + n(m-1)u^2\lv\nabla v\rv^2 ,
\end{align}
\end{subequations}
where $T_0=T-\frac{1}{n}\tr T\,g$ denotes the tracefree part of a symmetric $(0,2)$ tensor and all derivatives and traces are computed with respect to $g$.
\end{prop}

\begin{remark}

The condition that $(M^n,\hat g,\hat v^m\dvol_{\hat g})$ be quasi-Einstein is captured by~\eqref{eqn:qe_tf} and~\eqref{eqn:qe_lambda}, while the integrability condition~\eqref{eqn:kk_csc} is captured by~\eqref{eqn:qe_mu}.
\end{remark}

Note in particular that~\eqref{eqn:qe_conf} makes sense when we allow $u$ or $v$ to change signs.  From the standpoint of the tractor calculus, this is the more natural perspective to take (cf.\ \cite{Gover2008}).  While we shall always insist that an SMMS $(M^n,g,v^m\dvol_g)$ or an SCMS $(M^n,g,v^md\nu)$ be such that $v>0$, we shall otherwise not impose any sign restrictions on $u$ or $v$ (cf.\ Section~\ref{sec:formulation}).

Another important observation about~\eqref{eqn:qe_conf} is that the conditions~\eqref{eqn:qe_lambda} and~\eqref{eqn:qe_mu} imply that one can regard $\lambda$ and $\mu$ as the squared lengths of $u$ and $v$, respectively.  More precisely, if one sets $\hat u=cu$ and $\hat v=kv$ for constants $c,k>0$, then~\eqref{eqn:qe_conf} still holds provided one replaces $\lambda$ by $\hat\lambda=c^2\lambda$ and $\mu$ by $\hat\mu=k^2\mu$.  For this reason, when we consider fixing the measure $v^m\dvol_g$, we will frequently also fix the desired \charconstant\ $\mu$.  When we do this, we will say that we have fixed an SMMS $(M^n,g,v^m\dvol)$ with \charconstant\ $\mu$, so that the search for a function $u\in C^\infty(M)$ satisfying~\eqref{eqn:qe_conf} can be rephrased as the search for a quasi-Einstein scale as follows:

\begin{defn}
Let $(M^n,g,v^m\dvol_g)$ be an SMMS with \charconstant\ $\mu$.  Then $u\in C^\infty(M)$ is a \emph{quasi-Einstein scale} if the SMMS
\[ \left( M^n,u^{-2}g,u^{-m-n}v^m\dvol_g\right) \]
is a quasi-Einstein SMMS with \charconstant\ $\mu$ wherever it is defined.

The \emph{quasi-Einstein constant $\lambda$} of such a quasi-Einstein scale is the constant $\lambda$ determined by~\eqref{eqn:qe_conf}.
\end{defn}

Note in particular that this definition allows $u$ to change signs.  Similar terminology will be used for SCMS $(M^n,g,v^md\nu)$.

Finally, we observe that there is a ``duality'' between SMMS for certain pairs of dimensional parameter $m$.  This is the mechanism by which our approach to the study of quasi-Einstein metrics is seen to generalize the approach of He, Petersen and Wylie~\cite{HePetersenWylie2010,HePetersenWylie2011c}.

\begin{prop}[\cite{Case2010a}]
\label{prop:duality}
Let $(M^n,g,v^m\dvol_g)$ be an SMMS with \charconstant\ $\mu$ which admits a positive quasi-Einstein scale $u\in C^\infty(M)$ with quasi-Einstein constant $\lambda$.  Then the SMMS $(M^n,g,u^{2-m-n}\dvol_g)$ is an SMMS with \charconstant\ $\lambda$ such that $v$ is a quasi-Einstein scale with quasi-Einstein constant $\mu$.
\end{prop}

The proof hinges upon the observation that~\eqref{eqn:qe_conf} is invariant under the change of variables
\[ \left(u,v,m,\lambda,\mu\right) \mapsto \left(v,u,2-m-n,\mu,\lambda\right) . \]

Recast in terms of the warped products discussed in the context of Proposition~\ref{prop:kk}, this shows that given a Riemannian manifold $(M^n,g)$ and a fiber $(F^m,h)$, the problem of finding functions $u,v\in C^\infty(M)$ such that $\left(M^n\times F^m,\og=u^{-2}(g\oplus v^2h)\right)$ is Einstein is equivalent to either finding a quasi-Einstein scale for the SMMS $(M^n,g,v^m\dvol_g)$ or finding a quasi-Einstein scale for the SMMS $(M^n,g,u^{2-m-n}\dvol_g)$.  The only difference is in which function we are left to solve for --- $u$ or $v$, respectively --- and the meaning of the quasi-Einstein constant and the \charconstant\ for the metric $\og$.  In particular, this gives us the option to consider the SMMS $(M^n,g,1^{2-m-n}\dvol_g)$ with \charconstant\ $\lambda$ and to search for all quasi-Einstein scales $v$.  This is precisely the perspective of~\cite{HePetersenWylie2010,HePetersenWylie2011c}, but without our discussion of conformally invariant objects (and in particular, tractor bundles).

%% file: formulation.tex
\section{Tractor Formulation of Quasi-Einstein Metrics}
\label{sec:formulation}

Let us now turn to the problem of finding a formulation of quasi-Einstein metrics in terms of tractor bundles.  There are two approaches we can take to this problem, which have their own relative merits.

First, we can start with initial data a conformal manifold $(M^n,c)$ and ask whether, given a dimensional parameter $m$, there exist $u,v\in\mE[1]$ such that the system~\eqref{eqn:qe_conf} holds for some constants $\lambda$ and $\mu$.  Solutions to this problem will automatically give quasi-Einstein SMMS, and since this perspective does not impose any additional data than the dimensional parameter $m$, it is arguably the most natural perspective from the point of view of conformal geometry.  While we are not presently able to fully carry out this approach, we shall discuss some steps towards its resolution in Section~\ref{sec:formulation/new}.

Second, we can instead start with initial data an SCMS $(M^n,c,v^md\nu)$ with \charconstant\ $\mu$, and ask whether there exists a quasi-Einstein scale $u\in\mE[1]$.  From the perspective of SMMS, this is the most natural approach to the problem, as we are using the extra data of the measure in an important and essential way.  So long as we restrict $m\not\in\{-n,1-n,2-n\}$, we can completely carry out this approach.  In particular, we will describe in Section~\ref{sec:formulation/linear} a connection $\nabla^W$ on $\bT$ and a codimension one subbundle $\bT^W$ with the property that quasi-Einstein scales $u\in\mE[1]$ are in one-to-one correspondence with sections $I\in\mT^W$ which are parallel with respect to $\nabla^W$.

Unless otherwise stated, we assume for the rest of this article that all SCMS $(M^n,c,v^md\nu)$ are such that $m\not\in\{-n,1-n,2-n\}$.

\input{formulation/new}
\input{formulation/linear}

%% file: formulation/new.tex
\subsection{Tractor formulation}
\label{sec:formulation/new}

The purpose of this section is to show that the following definition is a generalization of the definition of a quasi-Einstein SCMS in the sense of Proposition~\ref{prop:conf_qe}.

\begin{defn}
\label{defn:tractor_qe}
Let $(M^n,c)$ be a conformal manifold with $n\geq 3$ and fix $m\in\bR$.  We say that $(I,J)\in\mT\oplus\mT$ is \emph{a quasi-Einstein pair} if
\begin{align}
\label{eqn:tractor_qe1} 0 & = (m+n-2)\lp X,J\rp\nabla I - m\lp X,I\rp\nabla J \\
\notag & \quad - \frac{1}{n}\big((m+n-2)\lp\nabla I,J\rp-m\lp\nabla J,I\rp\big)X \\
\label{eqn:tractor_qe2} 0 & = mn(m+n-2)\lp I,J\rp \big( \lp X,J\rp \nabla\lp X,I\rp - \lp X,I\rp\ \nabla\lp X,J\rp\big) \\
\notag & \quad - m(m+n-2)(2m+n-2)\lp X,I\rp \lp X,J\rp \nabla\lp I,J\rp \\
\notag & \quad + 2(m-1)(m+n-1)(m+n-2)\lp X,I\rp \lp X,J\rp \lp\nabla I,J\rp \\
\notag & \quad + 2m(m-1)(m+n-1)\lp X,I\rp \lp X,J\rp \lp\nabla J,I\rp .
\end{align}
\end{defn}

\begin{prop}
\label{prop:qe_implies_tractor}
Let $(M^n,c)$ be a conformal manifold and suppose that $u,v\in\mE[1]$ satisfy~\eqref{eqn:qe_conf}.  Then $I=\frac{1}{n}\bD u$ and $J=\frac{1}{n}\bD v$ are such that $(I,J)$ is a quasi-Einstein pair.
\end{prop}

\begin{proof}

The proof proceeds by rewriting the conditions~\eqref{eqn:qe_conf} in terms of $I$ and $J$.  For convenience, fix a scale $g\in c$ and write $I=(x,\nabla u,u)$ and $J=(y,\nabla v,v)$; i.e.\ the functions $x$ and $y$ are defined by
\[ x = -\frac{1}{n}\left(\Delta u + \J u\right), \qquad y = -\frac{1}{n}\left(\Delta v+\J v\right) . \]

First we isolate the condition~\eqref{eqn:qe_tf}.  Since $I$ and $J$ are scale tractors, the ``top'' component of~\eqref{eqn:tractor_qe1} automatically vanishes.  Writing~\eqref{eqn:qe_conf} in terms of the Schouten tensor yields
\begin{equation}
\label{eqn:qe_tf_p}
0 = (m+n-2)v\left(uP + \nabla^2 u + xg\right) - mu\left(vP+\nabla^2v+yg\right),
\end{equation}
establishing that the ``middle'' component of~\eqref{eqn:tractor_qe1} vanishes.  By the definitions of $x$ and $y$, it follows that
\begin{align*}
\delta\left(uP+\nabla^2 u+xg\right) & = (1-n)\big(dx-P(\nabla u)\big) \\
\delta\left(vP+\nabla^2v+yg\right) & = (1-n)\big(dy-P(\nabla v)\big) .
\end{align*}
Taking the divergence of~\eqref{eqn:qe_tf_p}, we then conclude that
\begin{align*}
0 & = (n-1)\left[(m+n-2)v(dx-P(\nabla u)) - mu(dy-P(\nabla v))\right] \\
& \quad -(m+n-2)(uP+\nabla^2 u+xg)(\nabla v) + m(vP+\nabla^2v+yg)(\nabla u) \\
& = n\left[(m+n-2)v(dx-P(\nabla u))-mu(dy-P(\nabla v))\right] \\
& \quad - (m+n-2)\lp\nabla I,J\rp + m\lp\nabla J,I\rp,
\end{align*}
thus establishing that the ``bottom'' component of~\eqref{eqn:tractor_qe1} vanishes.  Thus~\eqref{eqn:tractor_qe1} vanishes, as claimed.

To establish~\eqref{eqn:tractor_qe2}, we must also use the conditions~\eqref{eqn:qe_lambda} and~\eqref{eqn:qe_mu}.  To start, observe that these conditions are equivalently written in terms of $I$ and $J$ as
\begin{align}
\label{eqn:norm_u_tractor}
\lambda \lp X,J\rp^2 & = -(m+n-1)\lp X,J\rp^2 \lv I\rv^2 + m \lp X,I\rp \lp X,J\rp \lp I,J\rp \\
\label{eqn:norm_v_tractor}
\mu \lp X,I\rp^2 & = -(m+n-2)\lp X,I\rp \lp X,J\rp \lp I,J\rp + (m-1)\lp X,I\rp^2\lv J\rv^2 ,
\end{align}
respectively.  Differentiating~\eqref{eqn:norm_v_tractor} yields
\begin{align*}
0 & = -(m+n-2)\lp X,I\rp \lp X,J\rp \nabla\lp I,J\rp + 2(m-1)\lp X,I\rp^2 \lp\nabla J,J\rp \\
& \quad + (m+n-2)\lp I,J\rp \big(\lp X,J\rp\nabla\lp X,I\rp - \lp X,I\rp\nabla\lp X,J\rp\big) ,
\end{align*}
while~\eqref{eqn:tractor_qe1} implies that
\[ 0 = (m+n-2)(n-1)\lp X,J\rp \lp\nabla I,J\rp - mn\lp X,I\rp\lp\nabla J,J\rp + m\lp X,J\rp\lp\nabla J,I\rp \]
by taking the inner product with $J$.  Combining these equations yields~\eqref{eqn:tractor_qe2}, as desired.
\end{proof}

We conclude this section with two comments which may be of use in further understanding Definition~\ref{defn:tractor_qe}.

First, from the proof of Proposition~\ref{prop:qe_implies_tractor}, we see that~\eqref{eqn:qe_tf} alone implies the tractor formula~\eqref{eqn:tractor_qe1}.  Geometrically, \eqref{eqn:qe_tf} states that the SMMS
\[ (M^n,u^{-2}\bg,(v/u)^m\dvol_{u^{-2}\bg}) \]
has tracefree Bakry-\'Emery Ricci curvature.  This condition was considered by Catino~\cite{Catino2010}, who showed that under the additional assumption that the conformal class $c$ is locally conformally flat, such an SMMS necessarily splits locally as the warped product of an interval with a spaceform.

Second, it does not seem to be the case that~\eqref{eqn:tractor_qe1}, or even~\eqref{eqn:tractor_qe1} and~\eqref{eqn:tractor_qe2} together, imply that the tractors $I$ and $J$ are scale tractors; i.e.\ that $I=\frac{1}{n}\bD u$ and $J=\frac{1}{n}\bD v$ for some $u,v\in\mE[1]$.  However, we do not presently have an example which illustrates this behavior.

%% file: formulation/linear.tex
\subsection{The linear problem}
\label{sec:formulation/linear}

Let us now turn to the second approach to studying quasi-Einstein metrics using the tractor calculus, namely by fixing an SCMS $(M^n,c,v^md\nu)$ with \charconstant\ $\mu$.  The purpose of this section is to define the tractor subbundle $\bT^W$ and the operators $\nabla^W$ and $\bD^W$ appearing in Theorem~\ref{thm:equivalence}.

\begin{defn}
Let $(M^n,c,v^md\nu)$ be an SCMS with \charconstant\ $\mu$.  With $J=\frac{1}{n}\bD v$, the \emph{SCMS-scale tractor $\tilde J$} is defined by
\begin{equation}
\label{eqn:jtilde}
\tilde J = J + \frac{(m+2n-2)(\mu-(m-1)\lv J\rv^2)}{2(m+n-1)(m+n-2)\lp X,J\rp} X .
\end{equation}
Given a choice of scale $g\in c$, we will denote by $\tilde y$ the ``bottom'' component of $\tilde J$; i.e.\ $\tilde J=(\tilde y,\nabla v,v)$ for
\[ \tilde y = -\frac{2(n-1)v\Delta v+Rv^2-(m+2n-2)(\mu-(m-1)\lv\nabla v\rv^2)}{2(m+n-1)(m+n-2)v} . \]
\end{defn}

\begin{defn}
Let $(M^n,g,v^m\dvol_g)$ be an SMMS with \charconstant\ $\mu$.  The \emph{weighted Schouten tensor $P^W$} and the \emph{weighted Schouten scalar $\J^W$} are defined by
\begin{align*}
P^W & = \frac{1}{m+n-2}\left(\Ric_\phi^m - \J^Wg\right) \\
\J^W & = \frac{1}{2(m+n-1)}\left(R_\phi^m+m\mu v^{-2}\right) .
\end{align*}
\end{defn}

\begin{defn}
\label{defn:qe}
Let $(M^n,c,v^md\nu)$ be an SCMS with \charconstant\ $\mu$.  The \emph{$W$-tractor subbundle $\bT^W$} is the codimension one subbundle $\bT^W=\tilde J^\perp\subset\bT$.

The \emph{$W$-tractor connection $\nabla^W$} is defined given a choice of scale $g\in c$ by
\[ \nabla^W\begin{pmatrix}\sigma\\\omega\\\rho\end{pmatrix} = \begin{pmatrix}\nabla\sigma-\omega\\\nabla\omega+\sigma P^W+\rho g\\\nabla\rho - P^W(\omega)\end{pmatrix}. \]

A section $I\in\mT$ is \emph{$W$-parallel} if $\nabla^WI=0$.

The \emph{$W$-tractor-$D$ operator $\bD^W\colon\mE[w]\to\mT[w-1]$} is defined by
\begin{equation}
\label{eqn:w_d}
\bD^W\sigma = \begin{pmatrix} w(m+n+2w-2)\sigma\\(m+n+2w-2)\nabla\sigma\\-(\Delta_\phi\sigma+w\J^W\sigma)\end{pmatrix} .
\end{equation}
\end{defn}

Note that $\bT^W$ is clearly a well-defined codimension one subbundle of $\bT$, as $\tilde J$ is a well-defined nonvanishing section of $\bT$.  However, it is not clear that the $W$-tractor connection and the $W$-tractor-$D$ operator are well-defined tractor operators; i.e.\ it remains to check that they satisfy the transformation law~\eqref{eqn:tractor_law}.  One could verify this by direct computation.  However, we find it more interesting to reformulate $\nabla^W$ and $\bD^W$ in terms of $\nabla$ and $\bD$ in an invariant way, as will be carried out in Section~\ref{sec:bundle_properties}.

Assuming for the moment that $\nabla^W$ and $\bD^W$ are well-defined, we conclude this section with two simple observations.  First, the $W$-tractor connection preserves the tractor metric.

\begin{lem}
\label{lem:h_parallel}
Let $(M^n,c,v^md\nu)$ be an SCMS with \charconstant\ $\mu$.  The $W$-tractor connection $\nabla^W$ preserves the tractor metric $h$; i.e.\ $\nabla^Wh=0$.
\end{lem}

\begin{proof}

Choose a scale $g\in c$.  A direct computation shows that for all $I\in\mT$ and $x\in TM$,
\[ \frac{1}{2}\nabla_x\lv I\rv^2 = \sigma\nabla_x\rho + \rho\nabla_x\sigma + \lp\nabla_x\omega,\omega\rp = \lp\nabla_x^W I,I\rp . \qedhere \]
\end{proof}

Second, the tractor metric and the $W$-tractor-$D$ operator give a simple formula for the weighted scalar curvature of the SMMS determined by a choice of scale for an SCMS.

\begin{lem}
\label{lem:weighted_csc}
Let $(M^n,c,v^md\nu)$ be an SCMS with \charconstant\ $\mu$.  Fix $u\in\mE[1]$ and denote by $I$ the weighted scale tractor $I=\frac{1}{m+n}\bD^Wu$.  Given any choice of scale $g\in c$, it holds that
\[ \lv I\rv^2 = -\frac{(R_\phi^m+m\mu v^{-2})u^2 + 2(m+n-1)u\Delta_\phi u - (m+n)(m+n-1)\lv\nabla u\rv^2}{(m+n)(m+n-1)} . \]
In particular, the SMMS $(M^n,u^{-2}\bg,(u^{-1}v)^m\dvol_{u^{-2}\bg})$ satisfies~\eqref{eqn:kk_csc} if and only if $\lv I\rv^2=-\frac{\lambda}{m+n-1}$.
\end{lem}

\begin{proof}

This follows directly from the definitions of $\bD^W$ and $\J^W$.
\end{proof}

%% file: bundle_properties.tex
\section{Properties of the $W$-Tractor Connection}
\label{sec:bundle_properties}

In order to effectively make use of the $W$-tractor connection --- and in particular, to complete the proof of Theorem~\ref{thm:equivalence} --- we need to understand some of its basic properties.  While there is a nice algebraic framework into which these bundles fit (cf.\ Appendix~\ref{sec:algebra}; \cite{BCEG2006,Hammerl2008}), we will proceed in this task by direct computation.  As a potential benefit of this approach, we will arrive at more general forms of certain curvature formulae which appear in~\cite{Case2010b,HePetersenWylie2010} in more specialized problems involving quasi-Einstein SMMS.  In particular, we hope these computations both shed further insight into role of the curvature formulae appearing in~\cite{Case2010b,HePetersenWylie2010}, and also that they might foster an appreciation for the simplifications afforded by the tractor approach.

\input{bundle_properties/computations}
\input{bundle_properties/properties}
\input{bundle_properties/parallel}

%% file: bundle_properties/computations.tex
\subsection{Preliminary computations}
\label{sec:bundle_properties/computations}

In order to achieve the aforementioned goals, we first need to carry out some tedious computations involving $P^W$, $\J^W$, and certain of their derivatives and divergences.  For simplicity, throughout this section we fix an SMMS $(M^n,g,v^m\dvol)$ with \charconstant\ $\mu$, and define the function $y$ as the ``bottom component'' of the scale tractor $J=\frac{1}{n}\bD v$; i.e.\ $y=-\frac{1}{n}(\Delta v + \J v)$.  Note that, by~\eqref{eqn:jtilde}, this forces $y$ and $\tilde y$ to be related by
\[ \tilde y = y + \frac{(m+2n-2)(\mu-(m-1)\lv J\rv^2)}{2(m+n-1)(m+n-2)v} . \]

\begin{lem}
\label{lem:J}
\begin{equation}
\label{eqn:J^C}
\J^W = \J + mv^{-1}\left(y+\frac{\mu-(m-1)\lv J\rv^2}{2(m+n-1)v}\right) .
\end{equation}
\end{lem}

\begin{proof}

By the definition of $J\in\mT$, we have that
\begin{align*}
R_{\phi}^m+m\mu v^{-2} & = R - 2mv^{-1}\Delta v - m(m-1)v^{-2}|\nabla v|^2 + m\mu v^{-2} \\
& = \frac{m+n-1}{n-1}R + 2mnv^{-1}y - m(m-1)v^{-2}|\nabla v|^2 + m\mu v^{-2} \\
& = \frac{m+n-1}{n-1}R + 2m(m+n-1)v^{-1}y + mv^{-2}(\mu - (m-1)|J|^2) .
\end{align*}
The result then follows from the definitions of $\J$ and $\J^W$.
\end{proof}

\begin{lem}
\label{lem:P}
\begin{equation}
\label{eqn:P^C}
P^W = \frac{n-2}{m+n-2}P - \frac{mv^{-1}}{m+n-2}\left(\nabla^2 v + y\,g + \frac{\mu-(m-1)\lv J\rv^2}{2(m+n-1)v}g\right) .
\end{equation}
In particular,
\begin{subequations}
\label{eqn:trP^C}
\begin{align}
\label{eqn:trP^C_J^C} \tr P^W & = \J^W - mv^{-1}\tilde y \\
\label{eqn:trP^C_J} & = \J - \frac{mn(\mu-(m-1)\lv J\rv^2)}{2(m+n-1)(m+n-2)v^2} .
\end{align}
\end{subequations}
\end{lem}

\begin{proof}

Using~\eqref{eqn:J^C} and the definition of $P^W$, we compute that
\begin{align*}
P^W & = \frac{1}{m+n-2}\left(\Ric - mv^{-1}\nabla^2v - \left(\J + mv^{-1}y + \frac{m(\mu-(m-1)|J|^2)}{2(m+n-1)v^2}\right)g\right) \\
& = \frac{n-2}{m+n-2}P - \frac{mv^{-1}}{m+n-2}\left(\nabla^2 v + y\,g + \frac{\mu-(m-1)|J|^2}{2(m+n-1)v}g\right) .
\end{align*}

To achieve the second two claims, first observe that we may rewrite
\begin{equation}
\label{eqn:P^C-P}
P^W = P - \frac{mv^{-1}}{m+n-2}\left(vP + \nabla^2 v + y\,g + \frac{\mu-(m-1)|J|^2}{2(m+n-1)v}g\right) .
\end{equation}
By definition, $\tr(vP+\nabla^2v+y\,g)=0$, and so we see that
\[ \tr P^W = \J - \frac{mn(\mu-(m-1)|J|^2)}{2(m+n-1)(m+n-2)v^2}, \]
yielding~\eqref{eqn:trP^C_J}.  Using \eqref{eqn:J^C} again yields~\eqref{eqn:trP^C_J^C}.
\end{proof}

For the remaining lemmas, we require a few words on our conventions.  First, we will define the curvature $\mR\in\Lambda^2T^\ast M\otimes\End(E)$ of a connection $\nabla$ on a vector bundle $E$ by
\[ \mR(x,z)(s) = -(d^\nabla)^2s(x,z) := -\nabla_x\nabla_zs+\nabla_z\nabla_xs+\nabla_{[x,z]}s \]
for all $x,z\in TM$ and $s\in V$.  This is opposite the more common conventions in the literature (especially the tractor literature), but we will use it because it more naturally fits our form-based perspective for working with the adjoint tractor bundle (cf.\ \cite{Case2011o}).  In this perspective, the Kulkarni-Nomizu product is just the usual wedge product on $T^\ast M\otimes T^\ast M$,
\[ (h\wedge k)(x,y,z,w) = h(x,z)k(y,w) + h(y,w)k(x,z) - h(x,w)k(y,z) - h(y,z)k(x,w), \]
and the trace $\tr A\in T^\ast M\otimes T^\ast M$ of $A\in\Lambda^2T^\ast M\otimes\Lambda^2T^\ast M$ is given by
\[ \tr A(x,z) = \sum_{i=1}^n A(e^i,x,e^i,z) \]
for any orthonormal basis $\{e^i\}$ of $T_pM$, and likewise for other sections of $\Lambda^kT^\ast M\otimes\Lambda^lT^\ast M$.  With these conventions, the divergence $\delta T\in\Lambda^{k-1}T^\ast M\otimes E$ of an $E$-valued $k$-form $T$ is
\[ \delta T = \sum_{i=1}^n \nabla_{e^i}T(e^i,\dotsc) . \]

\begin{lem}
\label{lem:dP}
Let $(E,\nabla)$ be the tangent bundle $TM$ with the Levi-Civita connection and let $d\colon\Lambda^kT^\ast M\otimes E\to\Lambda^{k+1}T^\ast M\otimes E$ be the ``twisted'' exterior derivative.  Regarding $P,P^W\in\Lambda^1T^\ast M\otimes E$, it holds that
\begin{subequations}
\label{eqn:dP^C}
\begin{align}
\label{eqn:dP^C_rm} dP^W & = \frac{n-2}{m+n-2}dP - \frac{mv^{-1}}{m+n-2}\bigg(-\Rm(\nabla v) + dy\wedge g - v^{-1}dv\wedge\nabla^2 v \\
\notag & \qquad - v^{-1}dv\wedge y\,g + v\,d\big(\frac{\mu-(m-1)\lv J\rv^2}{2(m+n-1)v^2}\big)\wedge g\bigg) \\
\label{eqn:dP^C_weyl} & = \frac{n-2}{m+n-2}dP - \frac{mv^{-1}}{m+n-2}\bigg(-W(\nabla v) + (dy-P(\nabla v))\wedge g \\
\notag & \qquad - v^{-1}dv\wedge (vP+\nabla^2 v+y\,g) + v\,d(\frac{\mu-(m-1)\lv J\rv^2}{2(m+n-1)v^2})\wedge g\bigg) .
\end{align}
\end{subequations}

In particular,
\begin{equation}
\label{eqn:divP^C} \delta_\phi P^W = d\J^W + mv^{-2}\tilde y\,dv .
\end{equation}
\end{lem}

\begin{proof}

\eqref{eqn:dP^C_rm} follows immediately by taking the exterior derivative of~\eqref{eqn:P^C}.  To establish~\eqref{eqn:dP^C_weyl} from~\eqref{eqn:dP^C_rm}, one uses the Weyl decomposition
\begin{equation}
\label{eqn:weyl_decomposition}
\Rm = W + P \wedge g .
\end{equation}

To establish~\eqref{eqn:divP^C}, observe that by definition, the trace of $dP^W$ satisfies
\begin{equation}
\label{eqn:trdP^C}
\tr dP^W = \delta P^W - d\tr P^W .
\end{equation}
Using the fact that the Weyl tensor $W$ and the Cotton tensor $dP$ are traceless, it follows from~\eqref{eqn:dP^C_weyl} that
\begin{align*}
\tr dP^W & = \frac{mv^{-1}}{m+n-2}\bigg(v^{-1}(vP+\nabla^2 v+y\,g)(\nabla v) + (n-1)(dy-P(\nabla v)) \\
& \qquad + (n-1)v\,d\left(\frac{\mu-(m-1)|J|^2}{2(m+n-1)v^2}\right)\bigg) \\
& = -mv^{-1}P^W(\nabla v) - \frac{mv^{-1}}{m+n-2}\bigg((m-1)v^{-1}\nabla_{\nabla v}\nabla v + (m-1)v^{-1}y\,dv \\
& \qquad - (n-1)dy + \frac{m(\mu-(m-1)|J|^2)}{2(m+n-1)v^2}\,dv - (n-1)v\,d\left(\frac{\mu-(m-1)|J|^2}{2(m+n-1)v^2}\right)\bigg) \\
& = -mv^{-1}\left(P^W(\nabla v) - dy - \frac{m+2n-2}{2(m+n-1)(m+n-2)}d\left(\frac{\mu-(m-1)|J|^2}{v}\right)\right),
\end{align*}
where in the second line we have used~\eqref{eqn:P^C}.  Since $\delta_\phi=\delta+mv^{-1}\imath_{\nabla v}$, the result then follows from~\eqref{eqn:trP^C_J^C} and~\eqref{eqn:trdP^C}.
\end{proof}

\begin{lem}
\label{lem:trA}
Define the \emph{weighted Weyl curvature $A^W$} by $A^W=\Rm-P^W\wedge g$.  Then
\[ \tr A^W = mv^{-1}\left(vP^W+\nabla^2 v + \tilde yg\right) . \]
\end{lem}

\begin{proof}

Using~\eqref{eqn:weyl_decomposition}, we may write
\begin{equation}
\label{eqn:weyl2}
A^W = W-(P^W-P)\wedge g .
\end{equation}
Since $\tr W=0$ and $\tr (h\wedge g) = (n-2)h + \tr h\,g $ for all $h\in S^2T^\ast M$, it follows that
\begin{align*}
\tr A^W & = -(n-2)(P^W-P) - \tr (P^W-P)\,g \\
& = mv^{-1}\left(vP^W + \nabla^2 v + yg + \frac{(m+2n-2)(\mu-(m-1)|J|^2)}{2(m+n-1)(m+n-2)v}g\right),
\end{align*}
where the second line follows from~\eqref{eqn:P^C} and~\eqref{eqn:trP^C_J}.
\end{proof}

\begin{lem}
\label{lem:divA}
\begin{equation}
\label{eqn:divA^C}
\delta_\phi A^W = (m+n-3)dP^W - mv^{-2}dv\wedge \left(vP^W+\nabla^2v+\tilde yg\right) .
\end{equation}
\end{lem}

\begin{proof}

Recall that $\delta W = (n-3)dP$ (cf.\ \cite{Besse}).  Using~\eqref{eqn:weyl2}, we then see that
\begin{equation}
\label{eqn:div_weyl2}
\delta A^W = (n-2)dP - dP^W - \delta\left(P^W-P\right)\wedge g .
\end{equation}

Now, using~\eqref{eqn:dP^C} and the decomposition $\Rm=A^W+P^W\wedge g$, we see that
\begin{align*}
(n-2)dP - dP^W & = (m+n-3)dP^W + mv^{-1}\bigg(-A^W(\nabla v) + (dy-P^W(\nabla v))\wedge g \\
& \quad - v^{-1}dv\wedge(vP^W+\nabla^2v+yg) + vd\left(\frac{\mu-(m-1)\lv J\rv^2}{2(m+n-1)v^2}\right)\wedge g\bigg) .
\end{align*}
On the other hand, \eqref{eqn:trP^C_J} and~\eqref{eqn:trdP^C} imply that
\[ \delta(P^W-P) = -mv^{-1}\left(P^W(\nabla v)-d\tilde y + \frac{nv}{m+n-2}d\left(\frac{\mu-(m-1)\lv J\rv^2}{2(m+n-1)v^2}\right)\right) . \]
Inserting these two formulae into~\eqref{eqn:div_weyl2} yields the desired result.
\end{proof}

\begin{lem}
\label{tr:divdP}
As sections of $\Lambda^2T^\ast M$,
\[ \delta_\phi dP^W = -mv^{-2}dv\wedge \left(d\tilde y - P^W(\nabla v)\right) . \]
\end{lem}

\begin{remark}
To view $\delta_\phi dP^W$ as a section of $\Lambda^2T^\ast M$ means to consider the composition
\[ \Lambda^1T^\ast M \otimes \Lambda^1T^\ast M \xrightarrow{d} \Lambda^2T^\ast M\otimes \Lambda^1T^\ast M \xrightarrow{\delta_\phi} \Lambda^2T^\ast M . \]
\end{remark}

\begin{proof}

For convenience, we shall denote a contraction in two components by pairs of ``vectors'' $e$; e.g.\ $\Ric(x,z)=\Rm(e,x,e,z)$.  We compute at a point $p\in M$ by fixing $x,z\in T_pM$ and extending $x,z,e$ by parallel transport to a neighborhood of $p$.  It is then easy to see that
\begin{equation}
\label{eqn:intermediary}
\begin{split}
\delta_\phi dP^W(x,z) & = \nabla_e\nabla_x P^W(z,e) - \nabla_e\nabla_z P^W(x,e) + mv^{-1}dP^W(x,z,\nabla v) \\
& = \nabla_x\nabla_e P^W(z,e) - \nabla_z\nabla_e P^W(x,e) + \lp P^W(z),\Ric(x)\rp \\
& \quad - \lp P^W(x),\Ric(z)\rp + mv^{-1} dP^W(x,z,\nabla v) \\
& = \nabla_x\delta_\phi P^W(z) - \nabla_z\delta_\phi P^W(x) + mv^{-2}P^W(z,\nabla v)\lp \nabla v,x\rp \\
& \quad - mv^{-2}P^W(x,\nabla v)\lp \nabla v,z\rp ,
\end{split}
\end{equation}
where in the second line we have used the fact that $P^W$ is symmetric, whence $\Rm(e,x,P^W(e),z)$ is symmetric and in the third line we have used the facts
\begin{align*}
\Ric-mv^{-1}\nabla^2 v & = (m+n-2)P^W + \J^Wg \\
dP^W(x,z,\nabla v) &= \nabla_x\left(P^W(z,\nabla v)\right) - \nabla_z\left(P^W(x,\nabla v)\right) \\
& \quad - \lp P^W(z),\nabla_x\nabla v\rp + \lp P^W(x),\nabla_z\nabla v\rp .
\end{align*}
Using~\eqref{eqn:divP^C} and using the fact that the Hessian of a function is symmetric, the conclusion immediately follows from~\eqref{eqn:intermediary}.
\end{proof}

%% file: bundle_properties/properties.tex
\subsection{Properties}
\label{sec:bundle_properties/properties}

Let us now turn to establishing some of the basic, yet important, properties of $\nabla^W$ and $\bD^W$.  First, to verify that $\nabla^W$ and $\bD^W$ are well-defined operators, we give formulae for them in terms of the canonical tractor connection and the usual tractor-$D$ operator.

\begin{lem}
\label{lem:nabla_well_defined}
Let $(M^n,c,v^md\nu)$ be an SCMS with \charconstant\ $\mu$.  Then the $W$-tractor connection is given by
\[ \nabla^W = \nabla + \frac{mv^{-1}}{m+n-2}\left(\nabla J + \frac{\mu-(m-1)\lv J\rv^2}{2(m+n-1)v}\nabla X\right)\wedge X , \]
where we regard the second summand on the right hand side as a section of $T^\ast M\otimes\End(\bT)$ via the natural identification $\End(\bT)\cong\mA$.  In particular, $\nabla^W$ is well-defined as an operator $\nabla^W\colon\mT\to T^\ast M\otimes\mT$.
\end{lem}

\begin{proof}

For a fixed choice of scale $g\in c$, a simple computation shows that for any function $\psi\in\mE[0]$,
\begin{align*}
\nabla^W I - \nabla I & = \begin{pmatrix}0\\\sigma(P^W-P)\\-(P^W-P)(\omega)\end{pmatrix} \\
\big(\left(\nabla J + \psi\nabla X\right)\wedge X\big)(I) & = \begin{pmatrix}0\\-\sigma\left(vP+\nabla^2v+(y+\psi)g\right)\\\left(vP+\nabla^2v+(y+\psi)g\right)(\omega)\end{pmatrix} .
\end{align*}
The result then follows from~\eqref{eqn:P^C-P}.
\end{proof}

\begin{lem}
\label{lem:d_well_defined}
Let $(M^n,c,v^md\nu)$ be an SCMS with \charconstant\ $\mu$.  Then the $W$-tractor-$D$ operator $\bD^W\colon\mE[w]\to\mT[w-1]$ satisfies
\[ (n+2w-2)\bD^W\sigma - (m+n+2w-2)\bD\sigma = -mv^{-1}\lp\bD\sigma,J+\frac{\mu-(m-1)\lv J\rv^2}{2(m+n-1)v}X\rp\otimes X . \]
In particular, $\bD^W$ is a well-defined operator.
\end{lem}

\begin{proof}

The formula relating $\bD^W$ and $\bD$ follows from a direct computation using the definitions of $\bD^W$ and $\bD$ together with~\eqref{eqn:J^C}.  Since this formula writes $\bD^W$ in terms of well-defined tractor operators, this yields that $\bD^W$ is well-defined for $n+2w-2\not=0$.  That $\bD^W$ is well-defined when $n+2w-2=0$ is easily checked by verifying that the formula~\eqref{eqn:w_d} for $\bD^W$ transforms according to~\eqref{eqn:tractor_law} when one changes scales $g\mapsto e^{2s}g$.
\end{proof}

Second, we have the following computational lemma which is the essential ingredient in proving Theorem~\ref{thm:equivalence}.

\begin{lem}
\label{lem:algebra}
Let $(M^n,c,v^md\nu)$ be an SCMS with \charconstant\ $\mu$.  Let $I\in\mT$ and fix a scale $g\in c$ to denote $\nabla^WI=(\beta,T,\alpha)$.  Then
\[ \delta_\phi T = d\tr T - (m+n-1)\beta + mv^{-1}d\lp I,\tilde J\rp + (P^W-\J^W\,g)(\alpha) , \]
where the weighted divergence $\delta_\phi T = \delta T - T(\cdot,\nabla\phi)$ is taken in the second factor of $T\in T^\ast M\otimes T^\ast M$.
\end{lem}

\begin{proof}

By definition,
\begin{align*}
\delta_\phi T & = \sigma\delta_\phi P^W + P^W(\nabla\sigma) + \delta_\phi\nabla\omega + d\rho + mv^{-1}\rho\,dv \\
d\tr T & = \sigma\nabla\tr P^W + \tr P^W\,d\sigma + d\delta\omega + n\,d\rho .
\end{align*}
The only term which we do not yet have a useful equation for is $\delta_\phi\nabla\omega$.  To that end, we observe that the identity
\[ \delta\nabla\omega = d\delta\omega + \Ric(\omega) \]
together with the definition of $P^W$ implies that
\[ \delta_\phi\nabla\omega = d\delta\omega + (m+n-2)P^W(\omega) + \J^W\omega + mv^{-1}d\lp\omega,\nabla v\rp . \]

In order to get the desired statement, we first observe that combining~\eqref{eqn:trP^C_J^C} with~\eqref{eqn:divP^C} yields
\[ \delta_\phi P^W = d\tr P^W + mv^{-1}d\tilde y . \]
It then follows that
\begin{align*}
\delta_\phi T & = d\tr T + (m+n-1)P^W(\omega) + P^W(\nabla\sigma - \omega) - (n-1)d\rho + \J^W\omega - \tr P^W\,d\sigma \\
& \quad + mv^{-1}\left(d\lp\omega,\nabla v\rp + \sigma\,d\tilde y + \rho\,dv\right) \\
& = d\tr T - (m+n-1)\beta + (P^W - \J^W\,g)(\alpha) \\
& \quad + mv^{-1}d\left(\lp\omega,\nabla v\rp + \sigma\tilde y + v\rho\right) \\
& = d\tr T - (m+n-1)\beta + mv^{-1}d\lp I,\tilde J\rp + \left(P^W - \J^Wg\right)(\alpha) ,
\end{align*}
where the second line uses~\eqref{eqn:J^C}.
\end{proof}

As immediate consequences of Lemma~\ref{lem:algebra}, we have the following two simple but important facts.

\begin{cor}
\label{cor:prolongation}
Let $(M^n,c,v^md\nu)$ be an SCMS with \charconstant\ $\mu$.  Suppose that $I\in\mT^W$ is such that $\nabla^WI=\beta\otimes X$ for some one-form $\beta\in T^\ast M[-1]$.  Then $\beta\equiv0$.
\end{cor}

\begin{proof}

Fix a scale $g\in c$ and write $\nabla^WI=(\beta,T,\alpha)$.  By the assumption $I\in\mT^W$, we have that $\lp I,\tilde J\rp=0$.  By the assumption $\nabla^WI=\beta\otimes X$, we have that $T=0$ and $\alpha=0$.  That $\beta=0$ then follows immediately from Lemma~\ref{lem:algebra}.
\end{proof}

\begin{cor}
\label{cor:const_tilde_j}
Let $(M^n,c,v^md\nu)$ be a connected SCMS with \charconstant\ $\mu$.  Suppose that $I\in\mT$ is parallel with respect to $\nabla^W$.  Then $\lp I,\tilde J\rp$ is constant.
\end{cor}

\begin{proof}

Fix a scale $g\in c$ and write $\nabla^WI=(\beta,T,\alpha)$.  Since $\nabla^WI=0$, it follows immediately from Lemma~\ref{lem:algebra} that $d\lp I,\tilde J\rp=0$.
\end{proof}

The following basic lemma is necessary for the final corollary of Lemma~\ref{lem:algebra}.

\begin{lem}
\label{lem:traceT}
Let $(M^n,c,v^md\nu)$ be an SCMS with \charconstant\ $\mu$.  Let $I\in\mT$ and, fixing a scale $g\in c$, denote $I=(\rho,\omega,\sigma)$ and $\nabla^WI=(\beta,T,\alpha)$.  Then
\begin{equation}
\label{eqn:traceT}
\tr T + mv^{-1}\lp I,\tilde J\rp = (m+n)\rho + \delta_\phi\omega + \sigma\J^W .
\end{equation}
\end{lem}

\begin{proof}

This follows immediately from~\eqref{eqn:trP^C_J^C}.
\end{proof}

Combining Lemma~\ref{lem:algebra} and Lemma~\ref{lem:traceT}, we have the following key ingredient in the proof of Theorem~\ref{thm:equivalence}.

\begin{cor}
\label{cor:prolongation2}
Let $(M^n,c,v^md\nu)$ be an SCMS with \charconstant\ $\mu$.  Let $\sigma\in\mE[1]$ and suppose that $I=\frac{1}{m+n}\bD^W\sigma$ is such that $\nabla^WI=\beta\otimes X$ for some one-form $\beta\in T^\ast M[-1]$.  Then $\beta\equiv0$ and $I\in\mT^W$.
\end{cor}

\begin{proof}

Fix a metric $g\in c$ and write $I=(\rho,\nabla\sigma,\sigma)$ and $\nabla^WI=(\beta,T,\alpha)$.  By assumption, $\alpha=0$ and $T=0$.  By~\eqref{eqn:traceT} and the fact that $I=\frac{1}{m+n}\bD^W\sigma$, it follows that $\lp I,\tilde J\rp=0$.  The conclusion then follows from Corollary~\ref{cor:prolongation}.
\end{proof}

%% file: bundle_properties/parallel.tex
\subsection{Properties of $W$-parallel tractors}
\label{sec:bundle_properties/parallel}

We conclude this section by proving Theorem~\ref{thm:equivalence} and understanding the zero sets of $W$-parallel tractors $I\in\mT^W$.

\begin{proof}[Proof of Theorem~\ref{thm:equivalence}]

To begin, let $u\in\mE[1]$ and set $I=\frac{1}{m+n}\bD^Wu\in\mT$.  Fix a scale $g\in c$ and write $\nabla^WI=(\beta,T,0)$.  Using the definitions of $\nabla^W$ and the weighted Schouten tensor $P^W$, we see that
\begin{equation}
\label{eqn:scale_to_tractor}
\begin{split}
(m+n-2)T -(m+n-1)u^{-1}\lv I\rv^2g & = u\Ric_\phi^m + (m+n-2)\nabla^2u \\
& \quad + \left(\Delta_\phi u - (m+n-1)u^{-1}\lv\nabla u\rv^2\right)g \\
-(m+n)(m+n-1)\lv I\rv^2 & = \left(R_\phi^m+m\mu v^{-2}\right) u^2 + 2(m+n-1)u\Delta_\phi u \\
& \quad - (m+n)(m+n-1)\lv\nabla u\rv^2 .
\end{split}
\end{equation}
In particular, if $u$ defines a quasi-Einstein scale, the formulae given in Lemma~\ref{lem:conf_trans} imply that $\lv I\rv^2=-\frac{\lambda}{m+n-1}$ and $T=0$.  It then follows from Corollary~\ref{cor:prolongation2} that $I\in\mT^W$ and $\nabla^WI=0$.

Conversely, let $I\in\bT^W$ be such that $\nabla^WI=0$.  Fix a scale $g\in c$ and write $I=(\rho,\omega,u)$.  The vanishing of the top component of $\nabla^WI$ implies that $\omega=\nabla u$, while the vanishing of the middle component of $\nabla^W$ together with the assumption $I\in\mT^W$ imply, via Lemma~\ref{lem:traceT}, that $I=\frac{1}{m+n}\bD^Wu$.  It then follows from~\eqref{eqn:scale_to_tractor} that $u$ determines a quasi-Einstein scale with quasi-Einstein constant $\lambda=-(m+n-1)\lv I\rv^2$ and \charconstant\ $\mu$, as desired.
\end{proof}

In order to use Theorem~\ref{thm:equivalence}, one would like to know that, generically, $W$-parallel sections of $\bT^W$ are determined by positive functions.  One way to establish such a result is to use the analyticity of the quasi-Einstein equation (cf.\ \cite{HePetersenWylie2010}).  However, we can in fact prove more and classify the types of singular sets which can arise, generalizing the known behavior for almost Einstein metrics~\cite{Gover2008}, static metrics~\cite{Corvino2000}, and warped product Einstein metrics with boundary~\cite{HePetersenWylie2010}.

\begin{thm}
\label{thm:singularity_sets}
Let $(M^n,c,v^md\nu)$ be an SCMS with \charconstant\ $\mu$.  Let $I\in\mT^W$ be nonvanishing and suppose that $\nabla^W I=0$.  Denote by $\mS$ the singularity set of $I$, $\mS=\{p\in M\colon\lp X,I\rp_p=0\}$.
\begin{enumerate}
\item If $\lv I\rv^2<0$, then $\mS=\emptyset$.
\item If $\lv I\rv^2=0$, then either $\mS=\emptyset$ or $m=0$ and $\mS$ consists of isolated points.
\item If $\lv I\rv^2>0$, then either $\mS=\emptyset$ or $\mS$ consists of disjoint totally umbilic hypersurfaces.  Moreover, if $m\not=0$, these hypersurfaces are totally geodesic in the scale $v=1$.
\end{enumerate}
\end{thm}

\begin{proof}

Since $I$ is nonvanishing and parallel with respect to $\nabla^W$, $I$ cannot vanish at any point.  Also, by Theorem~\ref{thm:equivalence}, we know that $I=\frac{1}{m+n}\bD^W u$ for some $u\in\mE[1]$, with $\mS=u^{-1}(0)$.

Let $p\in\mS$, so that $\lv I(p)\rv^2=\lv\nabla u(p)\rv^2$.  As $\lv I\rv^2$ is constant, this shows that $\lv\nabla u\rv^2$ is constant along $\mS$.  Since $c$ is Riemannian, we see that if $\lv I\rv^2<0$, then no such $p$ can exist, and so $\mS=\emptyset$.  Suppose then that $\lv I\rv^2\geq 0$ and that $\mS\not=\emptyset$, and observe that $\lv I\rv^2=0$ if and only if $\nabla u(p)=0$.  In either case, since $\nabla^WI=0$, it also holds that
\begin{equation}
\label{eqn:umbilic}
0 = \nabla^2 u - \frac{1}{m+n}\Delta_\phi u\,g
\end{equation}
at $p$.

Now, if $\lv I\rv^2>0$, then $\lv\nabla u\rv^2=\lv I\rv^2>0$ along $\mS$.  Together with~\eqref{eqn:umbilic}, this implies that $\mS$ is a totally umbilic codimension one hypersurface, and moreover, if $m\not=0$, it is totally geodesic in the scale $v=1$.

On the other hand, if $\lv I\rv^2=0$, then $\nabla u(p)=0$.  Since $I(p)\not=0$, this implies that $\Delta_\phi u(p)\not=0$.  If $m\not=0$, this contradicts~\eqref{eqn:umbilic}.  Thus $m=0$, whence~\eqref{eqn:umbilic} implies that $\nabla^2u$ has a definite sign.  In particular, all zeroes of $u$ are isolated points.
\end{proof}

\begin{remark}
The situation where one has isolated singularities can arise for quasi-Einstein pairs $(I,J)$ with arbitrary $m\in[0,\infty)$, as can be seen by considering quasi-Einstein metrics of the form $(I,I)$ for $I$ any parallel tractor with isolated singularities; for an explicit example, see Appendix~\ref{sec:model}.  However, Theorem~\ref{thm:singularity_sets} implies that if $(I,J)$ is a quasi-Einstein pair and $p\in M$ is an isolated singularity for $I$ (resp.\ $J$), then necessarily $\lp X,J\rp_p$ (resp.\ $\lp X,I\rp_p$) must vanish.
\end{remark}

%% file: holonomy.tex
\section{Curvature, Holonomy, and Applications for Quasi-Einstein Metrics}
\label{sec:holonomy}

With Theorem~\ref{thm:equivalence} in hand, it is immediately clear that the vector space of quasi-Einstein scales on an SCMS with fixed \charconstant\ is finite-dimensional.  However much more can be said by considering the (reduced) holonomy group corresponding to $(\bT,\nabla^W)$.  This is because the Lie algebra of this group depends strongly on the curvature of $\nabla^W$, whence if an SCMS with \charconstant\ $\mu$ admits ``many'' linearly independent quasi-Einstein metrics, then it has suitably trivial curvature, and in particular belongs to a rather restrictive class of manifolds (cf.\ \cite{HePetersenWylie2011c}).

\input{holonomy/framework}
\input{holonomy/results}

%% file: holonomy/framework.tex
\subsection{Definitions and basic results}
\label{sec:holonomy/framework}

We begin by recalling the definition of the holonomy group of a vector bundle.

\begin{defn}
Let $E\to M$ be a vector bundle with connection $\nabla^E$, and fix a base point $p\in M$.  The \emph{holonomy group of $(E,\nabla^E)$ with respect to the base point $p$} is the Lie group
\[ \Hol(p,E,\nabla^E) = \left\{ \mP_\gamma^{\nabla^E} \colon \gamma \in \mL(M,p,p) \right\} \subset \GL(E_p), \]
where $\mL(M,p,q)$ denotes the space of piecewise smooth curves $\gamma\colon[0,1]\to M$ such that $\gamma(0)=p$ and $\gamma(1)=q$, and $\mP_\gamma^{\nabla^E}\colon E_p\to E_p$ is the parallel displacement along $\gamma$; for each $v\in E_p$, $\mP_{\gamma}^{\nabla^E}(v)=\psi_v(1)$, where $\psi_v\in\gamma^\ast E$ is the unique parallel section of $\gamma^\ast E$ satisfying $\psi_v(0)=v$.

The \emph{reduced holonomy group of $(E,\nabla^E)$ with respect to the base point $p$} is the Lie group
\[ \Hol_0(p,E,\nabla^E) = \left\{ \mP_\gamma^{\nabla^E} \colon \gamma \in \mL_0(M,p,p) \right\} \subset \GL(E_p), \]
where $\mL_0(M,p,p)$ denotes the space of contractible loops $\gamma\in\mL(M,p,p)$.
\end{defn}

Provided $M$ is connected, the (reduced) holonomy groups are all isomorphic.  Indeed,
\[ \Hol(p,E,\nabla^E) = \mP_{\gamma^{-1}}^{\nabla^E} \circ \Hol(q,E,\nabla^E) \circ \mP_\gamma^{\nabla^E}, \]
where $\gamma\in\mL(M,p,q)$, and a similar statement holds for $\Hol_0(p,E,\nabla^E)$.

Our interest is in the (reduced) holonomy groups corresponding to the bundle $(\bT,\nabla^W)$.  Indeed, we are only interested in the isomorphism type of these holonomy groups, and so we will study \emph{the} $W$-tractor holonomy group
\[ \Hol^W := \Hol(p,\bT,\nabla^W) , \]
where $p\in M$ is some fixed, but arbitrary, base point of $M$.  Similarly, we will denote by $\Hol_0^W$ the reduced $W$-tractor holonomy group.  Note also that, since $\nabla^W$ preserves the tractor metric, we can identify the fiber $\bT_p$ with $\bR^{n+1,1}$ and regard $\Hol^W\subset\SO(n+1,1)$.

As stated above, our goal is to use the (reduced) $W$-tractor holonomy group to study the space of quasi-Einstein scales on an SCMS with fixed \charconstant.  We can do this because of the following two facts.  First, by definition of the $W$-tractor holonomy group, there is a correspondence between quasi-Einstein scales and holonomy-invariant elements of $\bT_p^W$.

\begin{prop}
\label{prop:holonomy_invariant}
Let $(M^n,c,v^md\nu)$ be an SCMS with \charconstant\ $\mu$.  Then there is an isomorphism
\[ \big\{ u\in\mE[1]\colon u\mbox{ is a quasi-Einstein scale} \big\} \cong \big\{ I_p\in\bT_p^W \colon \Hol^W\cdot I_p = I_p \big\} . \]
\end{prop}

\begin{proof}

If $u\in\mE[1]$ is a quasi-Einstein scale, then $I=\frac{1}{m+n}\bD^Wu$ is parallel with respect to $\nabla^W$, and in particular $I_p$ is invariant under the action of $\Hol^W$.  Conversely, if $I_p$ is invariant under the action of $\Hol^W$, there is a uniquely determined tractor $I\in\mT^W$ such that $\nabla^WI=0$ which is determined at $p$ by $I_p$.  The result then follows from Theorem~\ref{thm:equivalence}.
\end{proof}

Second, the Ambrose-Singer Theorem (see~\cite{KobayashiNomizu1963}) yields a relationship between the Lie algebra of $\Hol_0^W$ and the curvature of $\nabla^W$.

\begin{thm}[Ambrose-Singer Theorem]
\label{thm:ambrose_singer}
Let $(M^n,c,v^md\nu)$ be an SCMS with \charconstant\ $\mu$.  Let $\mR^W=-(d^{\nabla^W})^2\in\Lambda^2T^\ast M\otimes\mA$ denote the curvature of the $W$-tractor connection.  Then
\[ \mathfrak{hol}(p,\bT,\nabla^W) = \mathrm{span}\left\{ P_{\gamma^{-1}}^{\nabla^W} \circ \mR^{W}(x,y) \circ \mP_\gamma^{\nabla^W} \colon x,y\in T_q M, \gamma\in\mL(M,p,q) \right\} , \]
where $\mathfrak{hol}(p,\bT,\nabla^W)$ is the Lie algebra of $\Hol_0(p,\bT,\nabla^W)$.
\end{thm}

In order to use the Ambrose-Singer Theorem, we need to know the curvature of $\nabla^W$.

\begin{prop}
\label{prop:curvature}
Let $(M^n,c,v^md\nu)$ be an SCMS with \charconstant\ $\mu$.  Then
\[ \mR^W(x,y) = \begin{pmatrix}0&0&0\\-dP^W(x,y)&A^W(x,y)&0\\0&dP^W(x,y)&0\end{pmatrix} .\]
\end{prop}

\begin{proof}

By definition, if $(\beta,T,\alpha)\in T^\ast M\otimes\mT$, then
\[ \nabla_x^W\begin{pmatrix}\alpha\\T\\\beta\end{pmatrix} = \begin{pmatrix} \nabla_x\alpha-T(x)\\\nabla_xT + \alpha\otimes P^W(x) + \beta\otimes x\\\nabla_x\beta-P^W(x,T(\cdot))\end{pmatrix} . \]
The result then follows by setting $\nabla^W I=(\beta,T,\alpha)$ and antisymmetrizing.
\end{proof}

\begin{remark}
There are three interesting observations about this result.

First, the term $A^W$ is the projecting part of $\mR^W$.  In particular, $A^W$ is conformally invariant.

Second, on an SMMS $(M^n,g,1^m\dvol_g)$ with \charconstant\ $\mu$, the tensor $A^W$ is not new.  Indeed, for $m>0$, it is precisely the tensor $A$ which appears in~\cite{Case2010b}, while in the dual case $m<2-n$, it is the tensor $Q$ appearing in~\cite{HePetersenWylie2010}.

Third, if $\nabla^WI=0$, then $\mR^W(I)=0$, which in the ``middle'' slot reads
\begin{equation}
\label{eqn:middle_curvature}
0 = \sigma\,dP^W - A^W(\omega) .
\end{equation}
If we take $I=\frac{1}{m+n}\bD^W u$, as per Theorem~\ref{thm:equivalence}, so that $\sigma=u$, $\omega=\nabla u$, and we compute in the scale $v=1$, then~\eqref{eqn:divA^C} and~\eqref{eqn:middle_curvature} together imply
\[ 0 = \delta A^W - (m+n-3)u^{-1}A^W(\nabla u) . \]
This is precisely the divergence-free condition for the tensor $A$ found in~\cite{Case2010b} and for $Q$ found in~\cite{HePetersenWylie2010}.
\end{remark}

%% file: holonomy/results.tex
\subsection{Applications}
\label{sec:holonomy/results}

Let us now use the holonomy groups $\Hol^W$ and $\Hol_0^W$ to study the space of quasi-Einstein scales on a given SCMS with \charconstant\ $\mu$.  As a first step, we need to understand the geometric implications of the flatness of the $W$-tractor connection.

\begin{lem}
\label{lem:weighted_flat}
Let $(M^n,c,v^md\nu)$ be an SCMS with \charconstant\ $\mu$.  The $W$-tractor connection $\nabla^W$ is flat if and only if one of the following conditions hold:
\begin{enumerate}
\item $m=0$ and $(M^n,c)$ is locally conformally flat,
\item $m=1$, $\mu=0$, and $(M^n,v^{-2}\bg)$ is a spaceform, or
\item $m\not\in\{0,1\}$ and $(M^n,v^{-2}\bg)$ is a spaceform with constant sectional curvature $-\frac{\mu}{m-1}$.
\end{enumerate}
\end{lem}

\begin{proof}

In general, decomposing $A^W$ according to the Weyl decomposition of algebraic curvature tensors
\[ \lp \mR \rp = \lp W \rp + \lp \Ric_0 \rp + \lp \id \rp , \]
we see that
\begin{equation}
\label{eqn:A^W_decomp}
A^W = W + \frac{mv^{-1}}{m+n-2}\left(vP+\nabla^2v\right)_0\wedge g + \left(\frac{m(\mu-(m-1)\lv J\rv^2)}{2(m+n-1)(m+n-2)v}\right)g\wedge g .
\end{equation}

If $\nabla^W$ is flat, we have that $A^W=0$, which in particular implies that the Weyl curvature vanishes.  Additionally, if $m\not=0$, the $\lp\Ric_0\rp$ component implies that $v^{-2}\bg$ is an Einstein metric and the $\lp\id\rp$ component implies that $\mu-(m-1)\lv J\rv^2=0$, implying the three stated conditions.

Conversely, if any of the three conditions of the lemma hold, then~\eqref{eqn:A^W_decomp} implies that $A^W=0$.  Moreover, \eqref{eqn:P^C-P} implies that $P=P^W$, whence $dP^W$ vanishes.  Hence, by Proposition~\ref{prop:curvature}, $\nabla^W$ is flat.
\end{proof}

The main application of Lemma~\ref{lem:weighted_flat} is to understand the rigid case of the following sharp bound on the dimension of the vector space of quasi-Einstein metrics on a conformal manifold (cf.\ \cite{HePetersenWylie2011c}).

\begin{thm}
\label{thm:finite}
Let $(M^n,c,v^md\nu)$ be an SCMS with \charconstant\ $\mu$ and denote by $\mQ^W$ the vector space of $W$-parallel sections of $\bT^W$.
\begin{enumerate}
\item $\dim\mQ^W\leq n+1$, and if equality holds, then $\nabla^W$ is flat.
\item If $\nabla^W$ is flat and $M^n$ is simply connected, then $\dim\mQ^W=n+1$.
\end{enumerate}
\end{thm}

\begin{proof}

Since $\dim\bT^W=n+1$, it follows immediately that $\dim\mQ^W\leq n+1$.  Suppose additionally that $\dim\mQ^W=n+1$.  Since $\dim\bT=n+2$, it follows that for any $I_1,I_2\in\mT$, there is an $I\in\mQ^W$ such that $I_1\wedge I_2=I\wedge I_3$ for some $I_3\in\mT$.  However, since $\nabla^WI=0$, we must have $\mR^W(I)=0$, whence $\mR^W$ annihilates $I_2\wedge I_3$.  Thus $\mR^W\equiv 0$, as claimed.

Conversely, suppose that $\nabla^W$ is a flat connection and $M$ is simply connected.  By Lemma~\ref{lem:weighted_flat}, this implies that $J^W=J=\frac{1}{n}\bD v$ and $\mu-(m-1)\lv\nabla J\rv^2=0$, whence $\nabla^W=\nabla$ and $\bT^W$ is a parallel subbundle of $\bT$.  By the Ambrose-Singer theorem, $\mathfrak{hol}^W$ is trivial, whence a basis for $\bT_p$ can be integrated to a parallel basis of $\bT$, and likewise for the parallel subbundle $\bT^W\subset\bT$.  In particular, $\dim\mQ^W=\dim\bT^W=n+1$.
\end{proof}

%% file: conclusion.tex
\section{Conclusion}
\label{sec:conclusion}

Let us conclude with some remarks on the wider applicability of our work.  In particular, we compare our results in Section~\ref{sec:holonomy} to recent results of He, Petersen and Wylie~\cite{HePetersenWylie2011c}, observe how our work might fit into the general problem of finding a global classification of locally conformally flat quasi-Einstein metrics, notice that the tractor formulation yields immediately formulae for some natural conformally covariant operators on SMMS, and observe that our tractor formulation is also related to the critical metrics recently studied by Miao and Tam~\cite{MiaoTam2008,MiaoTam2011}.

\begin{remark}
\label{rk:holonomy}
The general principal behind the study of holonomy is that reductions of the holonomy group correspond to additional geometric structures.  For example, Armstrong~\cite{Armstrong2007} has shown that there is an analogue of the De Rham decomposition theorem for the conformal holonomy group $\Hol(\bT,\nabla)$, which in particular allows him to conclude that manifolds admitting multiple almost Einstein metrics admit a sort of product decomposition.  Using completely different methods, He, Petersen and Wylie~\cite{HePetersenWylie2011c} have established a similar decomposition for SCMS with fixed \charconstant\ which admit multiple linearly independent quasi-Einstein scales.  A natural question is whether these results can be proven using holonomy methods for $(\bT,\nabla^W)$.

A related question arises in relation to the obstructions to the existence of quasi-Einstein metrics found in~\cite{Case2011o}.  There it is shown that if $A^W\colon TM\to\Lambda^2T^\ast M\otimes TM$ is injective, then one can construct an obstruction to the existence of quasi-Einstein metrics which depends polynomially on $A^W$ and $dP^W$.  In fact, this construction is really built from the curvature $\mR^W$ of $\nabla^W$, and the nondegeneracy assumption on $A^W$ can be viewed, via the Ambrose-Singer theorem, as an assumption on the holonomy algebra $\mathfrak{hol}^W$.  One might then hope to use the holonomy group $\Hol(\bT,\nabla^W)$ to understand what happens when $A^W$ is not injective, a program which has been partially carried out in the case of the normal tractor connection by Alt~\cite{Alt2006}.
\end{remark}

\begin{remark}
\label{rk:lcf}
A general question for smooth metric measure spaces, which is still open in full generality, is to give a complete classification of locally conformally flat quasi-Einstein metrics.  Of course, this question is trivial in the case $m=0$ --- they are all spaceforms.  More recently, the global classification in the case $m=\infty$ and $\lambda\geq0$ has been found~\cite{CaoChen2009}.  In the general case, however, only \emph{local} classifications are known: One can only conclude that locally the manifold splits as a warped product of an interval with a spaceform~\cite{CMMR2010,HePetersenWylie2010}.  Given the natural way that the assumption of local conformal flatness enters into the tractor calculus (as the flatness of the canonical tractor connection), it is natural to expect that the tractor calculus might lend some insights into the problem of finding the global classification.

An important obstacle in tackling this problem is to understand the examples found by B\"ohm~\cite{Bohm1999}, which are quasi-Einstein metrics on $S^n$ with positive quasi-Einstein constant and positive \charconstant\ (see Appendix~\ref{sec:model} for more details).  In particular, his examples are only constructed when $m+n\leq 9$, and it is not known if this is an essential assumption.
\end{remark}

\begin{remark}
\label{rk:operator}
An important class of operators in conformal geometry are conformally covariant operators; i.e.\ operators between conformal density bundles.  The most famous of these are the conformal Laplacian, the Paneitz operator, and the GJMS operators, which play an important role in studying the analytic and topological properties of manifolds; for further discussion and references, we refer the reader to~\cite{Chang2005}.  One can derive the existence of these operators using tractor techniques, which yields easily calculable formulae for low orders~\cite{GoverPeterson2003}.  In particular, this allows one to easily check that for $k=1,2$, the operators $L_{k,\phi}^m\colon\mE[-\frac{m+n-k}{2}]\to\mE[-\frac{m+n+k}{2}]$ defined by
\begin{align*}
L_{2,\phi}^m & = -\Delta_\phi + \frac{m+n-2}{2}\J^W \\
L_{4,\phi}^m & = \Delta_\phi^2 + \delta_\phi \circ (4P^W - (m+n-2)\J^W\id) \circ d \\
& \quad + \frac{m+n-4}{2}\left(-\Delta_\phi\J^W - 2\lv P^W\rv^2 - 2mv^{-2}\tilde y^2 + \frac{m+n}{2}(\J^W)^2\right)
\end{align*}
are the weighted analogues of the conformal Laplacian and the Paneitz operator, respectively.  An interesting open question is to study how these operators are related to their Riemannian counterparts; some results of this type for the weighted conformal Laplacian $L_{2,\phi}^m$ can be found in~\cite{Case2010b}.
\end{remark}

\begin{remark}
\label{rk:scalar}
In the case $(m,\mu)=(1,0)$, quasi-Einstein metrics are exactly static metrics in general relativity, which are often expressed in the form
\begin{equation}
\label{eqn:L}
0 = -v\Ric + \nabla^2 v - \Delta v\,g =: Lv .
\end{equation}
Here, $L$ is the formal adjoint of the linearization $DR$ of the scalar curvature (cf.\ \cite{Corvino2000}).  In the present language, \eqref{eqn:L} specifies that $(M^n,g,v^1\dvol)$ is a quasi-Einstein smooth metric measure space with \charconstant\ zero.

On the other hand, Miao and Tam~\cite{MiaoTam2008} have studied the critical points of the volume functional $V\colon \mM_\gamma^K\to\bR$, where $\mM_\gamma^K$ is the space of metrics $g$ on a compact Riemannian manifold $M$ with boundary $\Omega$ such that $g$ restricts to $\gamma$ on $\Omega$ and the scalar curvature of $g$ is $K$.  They found that $g$ is critical if and only if there is a $v\in C^\infty(M)$ such that $v>0$ on $M\setminus\Omega$, $v=0$ on $\Omega$, and
\begin{equation}
\label{eqn:miao_tam}
Lv = g .
\end{equation}

Working on the interior of $M$ where $\tilde J$ is defined, \eqref{eqn:miao_tam} has an interesting reinterpretation in terms of $(\mT,\nabla^W)$.  Namely, fix the SCMS $(M^n,c,v^1d\nu)$ with \charconstant\ $\mu=0$, where the density $v$ is such that it is trivialized to $v$ as in~\eqref{eqn:miao_tam} in the scale $g\in c$.  Suppose additionally that $I\in\mT$ is parallel with respect to $\nabla^W$ and that $u:=\lp X,I\rp>0$.  By Lemma~\ref{lem:algebra}, we know that $\lp I,\tilde J\rp$ is constant; let us assume that $\lp I, \tilde J\rp=\frac{1}{n-1}$.  Fix now the scale $u=1$.  Then the assumptions $\lp I,\tilde J\rp=\frac{1}{n-1}$ and $\nabla^W I=0$ imply that
\begin{align*}
\frac{1}{n-1} & = vx - \frac{1}{n}\left(\Delta v + \frac{R}{2(n-1)}v\right) \\
0 & = \Ric - v^{-1}\nabla^2 v - \frac{R-2v^{-1}\Delta v}{2n} g + (n-1)xg ,
\end{align*}
respectively, where we have written $I=(x,0,1)$.  Together, these imply that
\[ -\Ric + v^{-1}\nabla^2 v - v^{-1}\Delta v\,g = (n-1)\left(x - \frac{1}{n}\left(v^{-1}\Delta v + \frac{R}{2(n-1)}\right)\right) g; \]
i.e.\ \eqref{eqn:miao_tam} holds.

The key point we wish to make with this example is that the $W$-tractor connection is useful in its own right, and not simply when restricted to sections of $\bT^W$.  In particular, this example suggests that the framework for studying smooth metric measure spaces outlined in this paper should have applications beyond questions related to quasi-Einstein metrics.
\end{remark}

%% file: model.tex
\section{Quasi-Einstein Metrics on $S^n$}
\label{sec:model}

In this appendix we discuss some important quasi-Einstein metrics on $S^n$, beginning by considering its space of parallel tractors.  As is well-known, on $S^n$ the space $\mQ\subset\mT$ of parallel tractors is $(n+2)$-dimensional.  For an explicit description of $\mQ$, regard $S^n\subset\bR^{n+1}$ as the boundary of the ball $B(0,1)$ and let $\{x_1,\dotsc,x_{n+1}\}$ denote the standard Cartesian coordinates on $\bR^n$.  This fixes the standard metric $g_0$ on $S^n$ as the pullback of the Euclidean metric on $\bR^{n+1}$ via the inclusion map.  In this scale, a basis for $\mQ$ is
\[ \mB=\left\{I_0,I_1,\dotsc,I_{n+1}\right\} := \left\{ \begin{pmatrix} 1\\0\\-\frac{1}{2}\end{pmatrix},\begin{pmatrix}x_1\\\nabla^T x_1\\\frac{1}{2}x_1\end{pmatrix},\dotsc,\begin{pmatrix}x_{n+1}\\\nabla^T x_{n+1}\\\frac{1}{2}x_{n+1}\end{pmatrix} \right\} , \]
where $\nabla^T$ is the tangential part of the usual gradient in $\bR^n$; i.e.\ $\nabla^T$ is the gradient on $(S^n,g_0)$.  As is easy to check, $\mB$ is indeed an orthonormal basis for $\mQ$, with $\lv I_0\rv^2=-1$ and $\lv I_i\rv^2=1$ for $i=1,\dotsc,n+1$.

Geometrically, parallel tractors $I\in\mQ$ are equivalent to conformally flat Einstein metrics $g=\lp X,I\rp^{-2}g_0$.  For example, $I_0+I_1$ is equivalent to the Ricci flat metric $(1+x_1)^{-2}g_0$ on $S^n\setminus\{x_1=-1\}$, which is the pullback of the Euclidean metric on $\bR^n$ via the stereographic projection which maps $0\in\bR^n$ to the point $\{x_1=1\}$ and the boundary of the unit ball in $\bR^n$ to the equator $\{x_1=0\}$ in $S^n$.

In order to produce quasi-Einstein metrics according to Definition~\ref{defn:tractor_qe}, first observe that if $(I,J)\in\mQ\oplus\mQ\subset\mT\oplus\mT$, then automatically~\eqref{eqn:tractor_qe1} holds.  For such a pair, we may write
\[ I = \sum_{i=0}^{n+1} a_i I_i , \qquad J = \sum_{i=0}^{n+1} b_i I_i \]
for constants $a_0,\dotsc,a_{n+1},b_0,\dotsc,b_{n+1}$ and consider the function $\lambda$ defined by
\[ \lambda\lp X,J\rp^2 = -(m+n-1)\lp X,J\rp ^2 \lv I\rv^2 + m\lp X,I\rp \lp X,J\rp \lp I,J\rp . \]
As long as we assume $J\not\equiv0$, $\lambda$ will be a well-defined function, and we know from Section~\ref{sec:formulation/new} that~\eqref{eqn:tractor_qe2} holds if and only if $\lambda$ is a constant.  However, it is straightforward to check that this holds if and only if either $I$ is a constant multiple of $J$, or if $I$ and $J$ are linearly independent.  In particular, we see that the condition~\eqref{eqn:tractor_qe1} is strictly weaker than the quasi-Einstein condition.

\begin{remark}
\label{rk:bohm}
Not all quasi-Einstein pairs $(I,J)$ lie in the subspace $\mQ\oplus\mQ$.  Indeed, the examples of quasi-Einstein metrics found by B\"ohm~\cite{Bohm1999} are not in $\mQ\oplus\mQ$.  In regards to Remark~\ref{rk:lcf}, it is then natural to ask if we can classify all quasi-Einstein pairs $(I,J)$ on $S^n$.
\end{remark}

In order to have quasi-Einstein metrics in the sense of Definition~\ref{defn:qe}, it is necessary to fix either $(J,\mu,m)$ or $(I,\lambda,m)$, corresponding to fixing the SMMS
\begin{align*}
\left( S^n, g_0, v^m\dvol_g\right) & \quad \text{with \charconstant\ $\mu$}, \quad \text{or} \\
\left( S^n, g_0, u^{2-m-n}\dvol_g\right) & \quad \text{with \charconstant\ $\lambda$},
\end{align*}
respectively.  By Lemma~\ref{lem:weighted_flat}, if one wants the corresponding space $\mQ^W$ to achieve its maximal dimension, one must choose the \charconstant\ by $\mu=(m-1)\lv J\rv^2$ for $J$ parallel or $\lambda=-(m+n-1)\lv I\rv^2$ for $I$ parallel, respectively.  For example, the cases
\begin{align*}
(J,\mu,m) & = (I_0,-(m-1),m) \\
(I,\lambda,m) & = (I_0,m+n-1,m)
\end{align*}
are such that $\dim\mQ^W=n+1$, and here the elements of $\mQ^W$ are elliptic Gaussians (cf.\ \cite{Case2010a}).

%% file: algebra.tex
\section{The Algebraic Perspective}
\label{sec:algebra}

As alluded to in Section~\ref{sec:bundle_properties/computations}, there is a useful way to uniquely characterize the $W$-tractor connection by exploiting further the algebraic structure of the tractor bundles.  This can be done by using the language of BGG sequences to prolong the quasi-Einstein equation~\eqref{eqn:qe_conf} for either $(v,\mu)$ or $(u,\lambda)$ fixed.  We will summarize here the key algebraic facts, and refer the reader to an article of Hammerl~\cite{Hammerl2008}, both for a very readable account of the case of prolonging the conformally Einstein metric using the BGG machinery, as well as for a general outline for how this can be done for larger classes of conformally covariant operators (including the quasi-Einstein condition).

We begin by reviewing the algebraic structure underlying a conformal manifold, as used to realize conformal manifolds as examples of parabolic geometries.  Let $(M^n,c)$ be a conformal Riemannian manifold with $n\geq 3$, and let $\mA$ denote the adjoint tractor bundle as in Definition~\ref{defn:adjoint}.  Note in particular that, given a choice of scale $g\in c$, we may write
\begin{equation}
\label{eqn:grading}
\mA \ni B(x,(A,a),z) = \begin{pmatrix} -a&-x&0\\z&A&x^\flat\\0&-z^\sharp&a\end{pmatrix} ,
\end{equation}
where $(x,(A,a),z)\in TM \oplus (\Lambda^2T^\ast M\oplus C^\infty(M)) \oplus T^\ast M=:\mA_{-1}\oplus\mA_0\oplus\mA_1$, and $\flat\colon TM\to T^\ast M$ and $\sharp\colon T^\ast M\to TM$ are the usual ``musical'' isomorphisms determined by $g$.  Moreover, the algebraic bracket on $\mA$ is such that $\{\mA_i,\mA_j\}\subset\mA_{i+j}$, where we adopt the convention that $\mA_i=\{0\}$ for $|i|>1$.  Thus we see that each fiber $\mA_p$ is isomorphic to the $|1|$-graded Lie algebra
\[ \kg := \mathfrak{so}(n+1,1) = \bR^n \oplus \mathfrak{co}(n) \oplus (\bR^n)^\ast =: \kg_{-1}\oplus\kg_0\oplus\kg_1 . \]
Also, it is important to note that, while~\eqref{eqn:grading} does not give rise to a canonical grading of $\mA$, the fact that $\mA=TM\semiplus(\Lambda^2T^\ast M[-2]\oplus \mE[0])\semiplus T^\ast M$ implies that the adjoint tractor bundle can be canonically described as a filtered tractor bundle,
\[ \mA = \mA^{-1} \supset \mA^0 \supset \mA^1, \]
where $\mA^1=\mA_1$ and $\mA^0\cong\mA_0\oplus\mA_1$ for any choice of scale $g\in c$; in particular, each fiber $\mA_p$ is isomorphic to the filtered Lie algebra
\[ \kg = \kg^{-1} \supset \kg^0 \supset \kg^1, \]
defined in the same way.  Together with the normal tractor connection (or even the $W$-tractor connection), this shows that $(M^n,c)$ gives rise to a parabolic geometry modeled on $(\SO(n+1,1),P)$, where $P=\{g\in G\colon \Ad(g)(\kg^0)\subset\kg^0)\}$ (see~\cite{BaumJuhl2010,CapSlovak2009} for additional details).

On the standard tractor bundle, there is a natural family of connections coming from the natural Weyl structure on the conformal class $c$.  These are precisely those connections which differ from the normal tractor connection by choosing for each scale $g\in c$ a tensor $P_g^\prime$ which transforms the same way as the Schouten tensor; i.e.\ $P_{e^{2s}g}^\prime-P_{e^{2s}g}=P_g^\prime-P_g$ for all $s\in C^\infty(M)$.  Such connections can be decomposed according to the grading of $\kg$ as
\[ \nabla_Y^\prime I = Y\cdot I + \nabla^g_Y I + P^\prime(Y)\cdot I, \]
where we view $P^\prime\in\Hom(TM,T^\ast M)$ and the action $\cdot$ is induced from the action of $\mA$ on $\mT$ by the (possibly scale-dependent) inclusions $TM,T^\ast M\hookrightarrow\mA$.  The key fact here is that two such connections differ only in the highest homogeneity, and the difference is measured by a symmetric tensor.

Canonical descriptions of connections and differential operators in parabolic geometry are typically made using the Lie algebra cohomology groups $H^\ast(\kg_{1},V)$ for any finite-dimensional representation $V$ of $\kg$, or on the corresponding associated vector bundles.  For instance, in conformal geometry, one defines the spaces $H^\ast(M,\bT)$ using Kostant's codifferentials
\[ \partial_{k+1}^\ast \colon \Lambda^{k+1}T^\ast M\otimes\bT\to\Lambda^kT^\ast M\otimes\bT \]
defined by
\[ \partial_{k+1}^\ast(\xi_0\wedge\dotso\wedge\xi_k\otimes I) = \sum_{i=0}^k (-1)^{i+1} \xi_0\wedge\dotso\wedge\widehat{\xi_i}\wedge\dotso\wedge\xi_k\otimes(\xi_i\cdot I) . \]
It is easy to check that $\partial^\ast\circ\partial^\ast=0$, and so one can consider
\[ \mH_k:=H^k(M,\bT)=\ker\partial_k^\ast/\im\partial_{k+1}^\ast \]
as usual.

In this language, we can restate the computations from Section~\ref{sec:bundle_properties/computations} as follows.  First, Lemma~\ref{lem:P} implies that for the $W$-tractor connection, it holds that
\begin{equation}
\label{eqn:conn_c}
\partial_1^\ast(\nabla^W I) + mv^{-1}\lp I,\tilde J\rp X = \begin{pmatrix} 0\\\omega-\nabla\sigma\\\delta_\phi\omega+\sigma \J^W+(m+n)\rho\end{pmatrix}
\end{equation}
for all $I=(\rho,\omega,\sigma)\in\mT$.  Similarly, Lemma~\ref{lem:dP} and Lemma~\ref{lem:trA} imply that
\begin{equation}
\label{eqn:curv_c} 0 = \partial_2^\ast(\mR^WI) - mv^{-1}\lp\nabla^W\tilde J,I\rp X
\end{equation}
for all $I\in\mT$.

In the language of the BGG machinery, \eqref{eqn:conn_c} implies that the operator $u\mapsto\frac{1}{m+n}\bD^Wu$ is a splitting operator $L_0\colon\mH_0\to\ker\partial_0^\ast$ \emph{in the bundle $\bT^W$}.  This then implies that the quasi-Einstein condition~\eqref{eqn:qe_conf} is the corresponding BGG operator $D_0\colon\mH_0\to\mH_1$.  Moreover, \eqref{eqn:curv_c} then yields the isomorphism between $\ker D_0$ and $W$-parallel tractors, as guaranteed by Theorem~\ref{thm:equivalence} (cf.\ \cite{Hammerl2008,HSSS2010}).

\begin{remark}
Lemma~\ref{lem:trA} also implies that
\[ \partial^\ast\mR^W - mv^{-1}\nabla^W\tilde J\wedge X = 0 \]
which provides another way to canonically construct the $W$-tractor connection (cf.\ \cite{BaumJuhl2010,CapSlovak2009}).
\end{remark}